# Séries Gevrey de type arithmétique, I.
# Théorèmes de pureté et de dualité

Par Yves André

## Introduction

Ceci est le premier volet d'un article promouvant un point de vue arithmétique sur les fonctions spéciales classiques.

En feuilletant les traités de fonctions spéciales, on peut faire les trois observations suivantes. La première, vieille d'au moins un siècle, est que les séries explicites qui apparaissent comme développements de Taylor (ou bien asymptotiques) "classiques" sont des *séries Gevrey* d'ordre rationnel. Rappelons qu'une série entière $\sum_{n \geq 0} a_n z^n$ est dite Gevrey d'ordre $s$ si la série $\sum_{n \geq 0} \frac{a_n}{(n!)^s} z^n$ a un rayon de convergence non nul. La théorie, développée par Watson (1911), puis Ramis (dès 1978) et d'autres auteurs, s'est avérée féconde bien au-delà de l'étude des propriétés à l'infini des équations différentielles "classiques".

La seconde observation, d'où le présent article tire sa source, est qu'outre ces propriétés analytiques, les séries "classiques" — du moins celles à paramètres rationnels — possèdent aussi de remarquables propriétés arithmétiques, que nous englobons dans le concept de *série Gevrey de type arithmétique*:

Nous dirons qu'une série entière $\sum_{n \geq 0} a_n z^n$ est Gevrey d'ordre $s \in \mathbb{Q}$ de type arithmétique, si ses coefficients $a_n$ sont des nombres algébriques, et s'il existe une constante $C > 0$ telle que pour tout $n \in \mathbb{N}$, les conjugués du nombre algébrique $\frac{a_n}{(n!)^s}$ sont de module inférieur à $C^n$, et le dénominateur commun[1] à $a_0 = \frac{a_0}{(0!)^s}, \ldots, \frac{a_n}{(n!)^s}$ est inférieur à $C^n$.

C'est par exemple le cas de toute série hypergéométrique généralisée, confluente ou non, à paramètres rationnels.

Un peu plus généralement, nous définissons les *séries Nilsson-Gevrey d'ordre $s$ de type arithmétique* comme étant les sommes finies de la forme

$$\sum_{\alpha, k, l} \lambda_{\alpha, k, l} z^\alpha (\log^k z) \, y_{\alpha, k, l}(z)$$

---

[1] C'est-à-dire le plus petit entier naturel $d_n > 0$ tel que $d_n a_m/(m!)^s$ soit entier algébrique pour tout $m \leq n$.

où les $\lambda_{\alpha,k,l}$ sont des nombres complexes, $\alpha$ des nombres rationnels, $k, l$ des entiers naturels, et $y_{\alpha,k,l}(z)$ des séries Gevrey d'ordre $s$ de type arithmétique; pour $s$ fixé, on montre que ces "séries" forment une $\mathbb{C}[z]$-algèbre différentielle, notée NGA$\{z\}_s$.

La troisième observation est que la plupart des équations différentielles linéaires à coefficients polynômiaux qui émaillent les traités de fonctions spéciales sont de nature très simple; par exemple, on constate qu'ou bien elles sont fuchsiennes sur la sphère de Riemann, ou bien elles n'ont de singularités qu'à l'origine et à l'infini, l'une des deux étant régulière.

Nous nous proposons de montrer que la présence d'*une* série Nilsson-Gevrey *de type arithmétique* parmi les solutions suffit à rendre compte de cette dernière observation, en la précisant. Les principaux résultats sont les suivants (pureté).

THÉORÈME DE PURETÉ.    *Soit $y$ un élément de* NGA$\{z\}_s$ *vérifiant une équation différentielle linéaire $\Psi y = 0$ à coefficients dans $\mathbb{C}[z]$. On suppose $\Psi$ d'ordre minimal (en $\frac{d}{dz}$). Alors:*

i) *Si $s \leq 0$, $\Psi$ admet une base de solutions dans* NGA$\{z\}_s$;

ii) *Si $s > 0$, $\Psi$ admet une base de solutions de la forme* $\exp\left(\alpha_i z^{-\frac{1}{s}}\right) \cdot y_i$, *avec* $y_i \in$ NGA$\{z\}_s$, $\alpha_i \in \overline{\mathbb{Q}}$.

En termes plus vagues, toutes les séries qui interviennent dans les "solutions en 0" sont purement du même type. Comme les séries Gevrey d'ordre $\leq 0$ convergent, on voit en particulier que $0$ est une singularité régulière si $s \leq 0$. Dans le cas $s > 0$, on voit que $\frac{1}{s}$ est la seule pente non nulle du polygone de Newton de $\Psi$ à l'origine, et que les facteurs déterminants sont monômiaux.

Le cas $s = 0$ (cas fuchsien) est maintenant bien compris par la théorie des $G$-fonctions, cf. §3 ci-dessous; le résultat principal, dû en partie à Chudnovsky et en partie au présent auteur, exprime la permanence du "type arithmétique" par prolongement analytique (pureté en tout point algébrique):

THÉORÈME DE PERMANENCE ($s = 0$). *Soient $y$ et $\Psi$ comme ci-dessus, et supposons en outre que $s = 0$. Alors pour tout nombre complexe algébrique $a$, $\Psi$ admet une base de solutions dans* NGA$\{z - a\}_0$ (*resp. dans* NGA$\{\frac{1}{z}\}_0$).

On s'intéressera davantage ici au cas $s \neq 0$; c'est alors la prise de développement asymptotique qui remplace le prolongement analytique, et donne lieu à une "dualité" $0 \leftrightarrow \infty$:

THÉORÈME DE DUALITÉ $(s \neq 0)$. *Soient $y$ et $\Psi$ comme ci-dessus, et supposons en outre que $s \neq 0$. Alors $\Psi$ a au plus deux singularités non triviales:[2] l'origine et l'infini, dont l'une est régulière. De plus*:

i) *Si $s < 0$, $\Psi$ admet une base de solutions de la forme $\exp\left(\beta_i z^{-\frac{1}{s}}\right) \cdot \hat{y}_i$, avec $\hat{y}_i \in \mathrm{NGA}\left\{\frac{1}{z}\right\}_{-s}$, $\beta_i \in \overline{\mathbb{Q}}$. A fortiori, pour toute direction $\theta$ sauf un nombre fini (les directions de Stokes), il existe un nombre algébrique $\alpha_\theta$ tel que $\exp(\alpha_\theta z^{-\frac{1}{s}}) \cdot y$ admette à l'infini dans la direction $\theta$ un développement asymptotique dans $\mathrm{NGA}\left\{\frac{1}{z}\right\}_{-s}$;*

ii) *Si $s > 0$, $\Psi$ admet une base de solutions dans $\mathrm{NGA}\left\{\frac{1}{z}\right\}_{-s}$. A fortiori, pour toute direction $\theta$, $y$ est le développement asymptotique à l'origine dans la direction $\theta$ d'une solution[3] de $\Psi$ dans $\mathrm{NGA}\left\{\frac{1}{z}\right\}_{-s}$.*

On a donc aussi "pureté" à l'infini, mais pour l'ordre $-s$. Dans le cas $s < 0$, on voit que $-\frac{1}{s}$ est la seule pente non nulle du polygône de Newton de $\Psi$ à l'infini, et que les facteurs déterminants sont monômiaux.

Ces résultats peuvent être illustrés par l'exemple de l'intégrale oscillante d'Airy [OS]:

$$\mathrm{Ai}(z) = \frac{1}{\pi} \int_0^\infty \cos\left(zt + \frac{t^3}{3}\right) dt, \text{ solution de l'équation différentielle}$$
$$\text{d'ordre 2 (minimal) } \frac{d^2}{dz^2}\mathrm{Ai}(z) = z\mathrm{Ai}(z),$$

qui ne présente qu'une seule singularité (à l'infini), de pente $3/2$.

A l'aide du symbole de Pochhammer $(a)_n = a(a+1)\cdots(a+n-1)$, le développement de Taylor de $\mathrm{Ai}(z)$ en $0$ s'écrit

$$\frac{1}{\sqrt[3]{9}\Gamma\left(\frac{2}{3}\right)} \sum_{n \geq 0} \frac{1}{9^n \left(\frac{2}{3}\right)_n n!} z^{3n} - \frac{1}{\sqrt[3]{3}\Gamma\left(\frac{1}{3}\right)} \sum_{n \geq 0} \frac{1}{9^n \left(\frac{4}{3}\right)_n n!} z^{3n+1}.$$

Les deux séries qui apparaissent forment une base de solutions de l'équation d'Airy; elles sont Gevrey d'ordre $-2/3$ de type arithmétique (on peut utiliser le fait que pour $a, b$ rationnels, le dénominateur commun de $a/b, \ldots, (a)_n/(b)_n$ a une croissance au plus géométrique en $n$).

D'autre part, dans tout secteur d'ouverture $< 2\pi/3$ bissecté par le demi-axe réel positif, $\exp\left(\frac{2}{3}z^{\frac{3}{2}}\right) \cdot \mathrm{Ai}(z)$ admet le développement asymptotique à

---

[2] Une singularité est dite triviale si l'opérateur différentiel $y$ admet une base de solutions holomorphes.

[3] Si $\theta$ n'est pas singulière, la $\frac{1}{s}$-sommation de $y$ au sens de Ramis fournit une telle solution canonique. Si $\theta$ est une direction singulière, la $\frac{1}{s}$-sommation médiane y pourvoit.

l'infini

$$\frac{1}{2\sqrt{\pi}}\left(\frac{1}{z}\right)^{1/4}\sum_{n\geq 0}\left(\frac{3}{4}\right)^{2n}\frac{\left(\frac{1}{6}\right)_{2n}\left(\frac{5}{6}\right)_{2n}}{2n!}\left(\frac{1}{z}\right)^{3n}$$
$$-\frac{1}{2\sqrt{\pi}}\left(\frac{1}{z}\right)^{3/4}\sum_{n\geq 0}\left(\frac{3}{4}\right)^{2n+1}\frac{\left(\frac{1}{6}\right)_{2n+1}\left(\frac{5}{6}\right)_{2n+1}}{(2n+1)!}\left(\frac{1}{z}\right)^{3n+1}$$

et on constate que les deux séries en $\frac{1}{z}$ qui apparaissent sont Gevrey d'ordre $+2/3$ de type arithmétique (les directions d'angle $\pm\pi/3$ sont des directions de Stokes).

La preuve des théorèmes ci-dessus requiert des outils d'analyse $p$-adique *et* d'analyse complexe. Les techniques $p$-adiques servent à l'étude du cas $s = 0$ ($G$-opérateurs différentiels, §3); des techniques complexes interviennent dans l'étude détaillée (4.3) de ce que nous appelons $E$-opérateurs différentiels, c'est-à-dire des transformés de Fourier-Laplace des $G$-opérateurs différentiels (§§4, 5). L'importance de ces $E$-opérateurs apparaît dans le résultat suivant (§6):

*Soit $y(z)$ une série Nilsson-Gevrey d'ordre $s \neq 0$ de type arithmétique, vérifiant une équation différentielle à coefficients polynômiaux. Alors $y(z^{-s})$ est solution d'un $E$-opérateur différentiel.*

Ce volet se termine par une brève application des résultats ci-dessus à l'arithmétique des équations aux différences.

### Plan

1. Séries Gevrey de type arithmétique
2. Holonomie. $G$-fonctions et $E$-fonctions
3. G-opérateurs différentiels (compendium)
4. $E$-opérateurs différentiels
5. Autour de Laplace. Preuve de 4.3 et 4.6
6. Pureté et dualité
7. Séries de factorielles et opérateurs aux différences finies

## 1. Séries Gevrey de type arithmétique

1.1. Soit $\underline{a} = (a_0, a_1, a_2, \cdots)$ une suite de nombres algébriques. On considère la condition suivante:

(**G**): Il existe une constante $C > 0$ telle que pour tout $n$,

(**G**)$_{\text{conj}}$: Les conjugués de $a_n$ sont de module inférieur à $C^n$,

(**G**)$_{\text{dén}}$: Le dénominateur commun à $a_0, a_1, a_2, \cdots a_n$ est inférieur à $C^n$.

PROPOSITION 1.1.1.     *Soit $s = p/q$ un nombre rationnel écrit comme fraction irréductible, avec $q$ positif. Les conditions suivantes sont équivalentes*:

i) *la suite de terme général $\frac{a_n}{(n!)^s}$ vérifie* (**G**),

ii) *la suite de terme général $\frac{a_n}{(\left[\frac{n}{q}\right]!)^p}$ vérifie* (**G**).[4]

*Ces conditions sur $\underline{a}$ sont stables par addition, par produit de Cauchy, et par les transformations $a_n \to a_{n-1}$, $a_n \to (n+1)a_{n+1}$, $a_n \to \frac{1}{n}a_{n-1}$.*

*Preuve.* En ce qui concerne (**G**)$_{\text{conj}}$, l'équivalence de i) et ii) vient de ce que, d'après Stirling, $\log([\frac{n}{q}]!)^p \approx \log(n!)^s - ns\log q + O(\log n)$. Notons que i) équivaut à dire que la suite de terme général $\frac{(a_n)^q}{(n!)^p}$ vérifie (**G**). En ce qui concerne (**G**)$_{\text{dén}}$, l'équivalence de i) et ii) se ramène donc à comparer le dénominateur commun de $(a_0)^q, \ldots, \frac{(a_n)^q}{(n!)^p}$, d'une part, et celui de $(a_0)^q, \ldots, \frac{(a_n)^q}{(\left[\frac{n}{q}\right]!)^{pq}}$ d'autre part. Or $([\frac{n}{q}]!)^q$ divise $(q[\frac{n}{q}])!$, qui divise lui-même $n!$, et il s'agit de montrer que le plus petit commun multiple de $1, 1, \ldots, \frac{n!}{([\frac{n}{q}]!)^q}$ a une croissance au plus exponentielle en $n$. On observe que seuls des nombres premiers $p \leq n$ divisent $\frac{m!}{([\frac{m}{q}]!)^q}$ pour $m \in \{1, 2, \ldots, n\}$, l'exposant étant

$$\left[\frac{m}{p}\right] - q\left[\frac{m}{pq}\right] + \left[\frac{m}{p^2}\right] - q\left[\frac{m}{p^2q}\right] + \cdots \leq q\frac{\log m}{\log p}.$$

On a donc

$$\operatorname*{ppcm}_{m \leq n} \frac{m!}{([\frac{m}{q}]!)^q} \leq (n^q)^{\text{Card}\{p \leq n\}},$$

et le théorème de Chebishev Card$\{p \leq n\} = O(n/\log n)$ suffit pour conclure.

Traitons ensuite la question des dénominateurs quant à la stabilité par produit de Cauchy: Si $s \geq 0$, on écrit, pour $n$ fixé,

$$\operatorname*{dén}_{m \leq n} \sum_{1 \leq m} \frac{a_1 b_{m-1}}{m!^s} = \operatorname*{dén}_{m \leq n} \sum_{1 \leq m} \binom{m}{1}^s \frac{a_1}{1!^s}\frac{b_{m-1}}{(m-1)!^s} \leq \operatorname*{dén}_{1 \leq m \leq n} \left(\frac{a_1}{1!^s} \cdot \frac{b_{m-1}}{(m-1)!^s}\right);$$

Si $s < 0$, on écrit

$$\operatorname*{dén}_{m \leq n} \sum_{1 \leq m} \left(\left[\frac{m}{q}\right]!\right)^{-p} a_1 b_{m-1} \leq \operatorname*{dén}_{1 \leq m \leq n} \left(\left(\left[\frac{m}{q}\right]!\right)^{-p} a_1 b_{m-1}\right)$$

$$\leq \operatorname*{dén}_{1 \leq m \leq n} \left(\left(\left[\frac{1}{q}\right]!\right)^{-p} a_1 \cdot \left(\left[\frac{m-1}{q}\right]!\right)^{-p} b_{m-1}\right);$$

dans les deux cas, on obtient la croissance exponentielle souhaitée.

---

[4] La constante $C$ impliquée dans (**G**) n'étant pas nécessairement la même dans i) et dans ii).

Pour la question des dénominateurs quant à la stabilité par transformation $a_n \to \frac{1}{n} a_{n-1}$, on invoque de nouveau Chebishev sous la forme

$$\log \mathrm{ppcm}(1, 2, \ldots, n) = O(n).$$

Les autres vérifications relatives à la seconde assertion de la proposition sont immédiates.

Nous utiliserons librement le lemme suivant familier en théorie des $G$-fonctions depuis Siegel.

LEMME 1.1.2.    *Soient $a, b$ deux nombres rationnels, non entiers négatifs. Alors le dénominateur commun de $a/b, \ldots, (a)_n/(b)_n$ a une croissance au plus géométrique en $n$.*

1.2. Lorsque les conditions de la proposition 1.1.1. sont réalisées, nous dirons que la série entière $f = f_{\underline{a}} = \sum_{n \geq 0} a_n z^n$ est *Gevrey d'ordre $s$ de type arithmétique* (ou, plus brièvement, que $f$ est une série Gevrey arithmétique d'ordre $s$).

La seconde assertion de 1.1.1 montre que les séries Gevrey d'ordre $s$ de type arithmétique forment une sous-$\mathbb{Q}[z]$-algèbre de $\overline{\mathbb{Q}}[[z]]$ stable par différentiation et intégration. Nous noterons $\overline{\mathbb{Q}}\{z\}_s^A$ cette algèbre différentielle, et $K\{z\}_s^A$ la sous-algèbre différentielle formée des séries dont tous les coefficients $a_n$ appartiennent à un corps de nombres donné $K$.

*Remarques.* a) La suite de terme général $\frac{a_n}{(n!)^s}$ vérifie $(\mathbf{G})_{\mathrm{conj}}$ si et seulement si le rayon de convergence des séries $\sum_{n \geq 0} \iota(a_n) z^n$ (pour tout plongement $\iota : \overline{\mathbb{Q}} \to \mathbb{C}$) est minoré par une constante $> 0$ indépendante de $\iota$. En particulier, chaque $\iota$ induit un plongement $\overline{\mathbb{Q}}\{z\}_s^A \to \mathbb{C}\{z\}_s$ dans l'algèbre différentielle des séries Gevrey d'ordre s usuelles.

b) Si $K$ est un corps de nombres plongé dans $\mathbb{C}$, l'inclusion de $K\{z\}_s^A$ dans $\mathbb{C}\{z\}_s$ jouit de la propriété suivante: *Tout élément de $K\{z\}_s^A$ est soit un polynôme, soit une série Gevrey d'ordre exactement $s$* (i.e. n'est Gevrey d'ordre $r$ pour aucun $r < s$). En particulier $K\{z\}_s^A \cap K\{z\}_{s'}^A = K[z]$ si $s \neq s'$.

Une formulation équivalente est que si une suite de terme général $a_n \in K$ vérifie $(\mathbf{G})$ et a une infinité de termes non nuls, la série $\sum_{n \geq 0} a_n z^N$ a un rayon de convergence fini. De fait, la norme d'un terme $a_n$ non nul est un nombre rationnel de dénominateur et numérateur bornés par $C^{n[K:\mathbb{Q}]}$ par hypothèse, ce qui entraîne bien que $\limsup |a_n|^{\frac{1}{n}} > 0$.

c) Rappelons que pour $s < 0$, les éléments de $\mathbb{C}\{z\}_s$ s'identifient aux fonctions entières $F$ d'ordre $-\frac{1}{s}$, i.e. pour lesquelles il existe une constante $B > 0$ telle que $|F(z)| = O(e^{B|z|^{-\frac{1}{s}}})$.

1.3. *Changement d'ordre.* L'étude des séries Gevrey de type arithmétique se ramène dans une large mesure à celle des séries d'ordre $s = -1, 0$ ou $+1$. Voici comment.

Soient $u, v$ deux entiers naturels, $u \neq 0$, et introduisons les transformations de séries formelles suivantes, familières en analyse ultramétrique:

$$\varphi_u \left( \sum_{n \geq 0} a_n z^n \right) = \sum_{n \geq 0} a_n z^{un}, \quad \Psi_{v/u} \left( \sum_{n \geq 0} a_n z^n \right) = \sum_{n \geq 0} a_{un+v} z^n.$$

PROPOSITION 1.3.1.    *Soit $s$ un nombre rationnel. Les conditions suivantes sont équivalentes*:

i) $f \in \overline{\mathbb{Q}}\{z\}_s^A$,

ii) $\varphi_u f \in \overline{\mathbb{Q}}\{z\}_{s/u}^A$,

iii) *pour tout* $v \in \{0, \ldots, u-1\}$, $\Psi_{v/u} f \in \overline{\mathbb{Q}}\{z\}_{su}^A$.

Seules les questions relatives aux dénominateurs sont non triviales. L'équivalence de i) et ii) découle alors du fait vu en 1.1.1 le plus petit commun multiple de $1, 1, \ldots, \frac{nu!}{n!^u}$ a une croissance au plus exponentielle en $n$. Pour l'équivalence de i) et iii), on applique i) $\Leftrightarrow$ ii) à la série $\Psi_{v/u} f$ (au lieu de $f$), en remarquant que

$$z^v \varphi_u \Psi_{v/u} f = \sum_{n \geq 0} a_{un+v} z^{un+v},$$

de sorte que $f = \sum_{v=0}^{u-1} z^v \varphi_u \Psi_{v/u} f$. Les détails sont laissés au lecteur.

COROLLAIRE 1.3.2.    *Soit $s = p/q$ un nombre rationnel non nul écrit comme fraction irréductible, (avec $q$ positif), et soit $\varepsilon(s)$ le signe de $s$. Alors on a*

$$\sum_{n \geq 0} a_n z^n \in \overline{\mathbb{Q}}\{z\}_s^A$$

*si et seulement si pour tout* $v \in \{0, \ldots, q-1\}$,

$$\sum_{n \geq 0} a_{qn+v} z^{|p|n} \in \overline{\mathbb{Q}}\{z\}_{\varepsilon(s)}^A.$$

*Remarque.* Le produit de Hadamard d'une série Gevrey d'ordre $r$ de type arithmétique par une série Gevrey d'ordre $s$ de type arithmétique est une série Gevrey d'ordre $r + s$ de type arithmétique. En revanche le produit de Cauchy de séries Gevrey de type arithmétique d'ordres différents n'est pas en général de type arithmétique.

1.4. Un plongement $\overline{\mathbb{Q}} \to \mathbb{C}$ étant fixé, nous définissons (comme dans l'introduction) les séries *Nilsson-Gevrey d'ordre $s$ de type arithmétique* comme étant les sommes finies de la forme

$$y(z) = \sum_{\alpha, k} \lambda_{\alpha, k, l} z^{\alpha} \log^k z \, y_{\alpha, k, l}(z),$$

avec

$$\lambda_{\alpha, k, l} \in \mathbb{C}, \ \ \alpha \in \mathbb{Q}, \ \ k, l \in \mathbb{N}, \ \text{et } y_{\alpha, k, l}(z) \in \overline{\mathbb{Q}}\{z\}_s^A.$$

On peut former l'algèbre différentielle compositum de $\mathbb{C}[\log z]$ et $\overline{\overline{\mathbb{Q}[[z]]}}$ (séries de Puiseux à coefficients algébriques). Il est alors clair que, pour $s$ fixé, les séries Nilsson-Gevrey d'ordre $s$ de type arithmétique forment une sous-$\mathbb{C}[z]$-algèbre différentielle de ce compositum. Nous la noterons $\mathrm{NGA}\{z\}_s$ .

PROPOSITION 1.4.1. *Soit $s$ un nombre rationnel non nul, et notons $\varepsilon(s)$ son signe. Les conditions suivantes sont équivalentes*:

i) $y(z) \in \mathrm{NGA}\{z\}_s$,

ii) $y(z^{|s|}) \in \mathrm{NGA}\{z\}_{\varepsilon(s)}$.

Cela découle directement du corollaire précédent.

## 2. Holonomie. $G$-fonctions et $E$-fonctions

2.1. Soit $\underline{a} = (a_0, a_1, a_2, \ldots)$ une suite de nombres algébriques. On considère la condition suivante:

(**H**): La suite $\underline{a}$ vérifie une relation de récurrence de la forme

$$P_0(n)a_n + P_1(n+1)a_{n+1} + \cdots + P_\mu(n+\mu)a_{n+\mu} = 0$$

où les $P_0, \ldots, P_\mu$ désignent des polynômes à coefficients algébriques.

On remarque que les termes d'une suite vérifiant (**H**) engendrent une extension finie de $\mathbb{Q}$; on notera souvent $K$ le corps de nombres engendré par ces termes et par les coefficients des $P_j$.

Il est bien connu que la condition (**H**) équivaut à dire que $f$ satisfait à une équation différentielle linéaire à coefficients polynômiaux (on dit aussi que $f$ est *holonome*).

PROPOSITION 2.1.1. *Soit $s = p/q$ un nombre rationnel écrit comme fraction irréductible, avec $q$ positif. Les conditions suivantes sont équivalentes*:

i) *La suite $\underline{a}$ vérifie* (**H**),

ii) *La suite de terme général $\frac{a_n}{(\lfloor \frac{n}{q} \rfloor!)^p}$ vérifie* (**H**).[5]

---

[5] Le rang $\mu$ de la récurrence impliquée dans (**H**) n'étant pas nécessairement le même dans i) et dans ii).

*Ces conditions sur $\underline{a}$ sont stables par addition, par produit de Cauchy, par produit terme à terme, et par les transformations $a_n \to a_{n-1}$, $a_n \to (n+1)a_{n+1}$, $a_n \to a_{n-1}$; de plus, si $u,v$ sont deux entiers naturels, $u \neq 0$, ces conditions sont aussi stables par les transformations $a_n \to (a_{n/u}$ si $u$ divise $n$, $0$ sinon), $a_n \to a_{un+v}$.*

*Preuve.* La seconde assertion (relativement à la condition i)) signifie respectivement que l'holonomie d'une série est préservée par somme, produit, produit de Hadamard (cf. [A, I.3]), multiplication par $z$, dérivation, intégration, et par les opérations $\varphi_u$ et $\Psi_{v/u}$, ce qui est bien connu (pour $\Psi_{v/u}$, on peut utiliser la formule $\Psi_{v/u}f = \sum_{\varepsilon^u = 1}(z^{-v}f)(\varepsilon z)$).

Déduisons-en la première assertion. Posons $g = \sum_{n \geq 0} \frac{a_n}{\left[\frac{n}{q}\right]!} z^n$. Alors $\Psi_{v/q}q = \sum_{n \geq 0} \frac{a_{qn+v}}{n!^p} z^n$, $\Psi_{v/q}f = \sum_{n \geq 0} a_{qn+v}z^n$, et en appliquant les formules $f = \sum_{v=0}^{q-1} z^v \varphi_q \Psi_{v/q}f$ et $g = \sum_{v=0}^{q-1} z^v \varphi_q \Psi_{v/q}g$, on se ramène (en itérant $|p|$ fois) à montrer que $f = \sum_{n \geq 0} a_n z^n$ est holonome si et seulement $F = \sum_{n \geq 0} \frac{a_n}{n!}z^n$ l'est. Cela résulte aisément du tableau de correspondances suivant:

$$(2.1.2) \qquad f(z) \leftrightarrow F(z),$$

$$(2.1.3) \qquad z\frac{d}{dz}f(z) \leftrightarrow z\frac{d}{dz}F(z),$$

$$(2.1.4) \qquad \left(z^2\frac{d}{dz} + z\right)f(z) \leftrightarrow zF(z),$$

$$(2.1.5) \qquad zf(z) \leftrightarrow \int_0^z F(z).$$

2.2. Soit $\underline{a}$ une suite de nombres algébriques vérifiant (**G**) et (**H**). On dit alors que la série

$$f = f_{\underline{a}} = \sum_{n \geq 0} a_n z^n$$

est *de type G*, ou est une *G-fonction*, et que la série

$$F = F_{\underline{a}} = \sum_{n \geq 0} \frac{a_n}{n!} z^n$$

est *de type E*, ou est une *E-fonction*.

Il y a lieu de considérer aussi la série

$$\mathfrak{f} = \mathfrak{f}_{\underline{a}} = \sum_{n \geq 0} n! a_n z^n.$$

Nous donnerons à une telle série le nom de série *de type Э*, ou *Э-fonction*[6] (elle diverge aux places archimédiennes mais définit des fonctions analytiques en 0 pour toute place ultramétrique, dont les propriétés diophantiennes ont été étudiées par V. Chirskii). Un exemple typique est la série d'Euler

$$\sum_{n \geq 0} (-1)^n n! z^{n+1}$$

qui est le développement asymptotique à l'origine de la fonction $\int_0^\infty \frac{e^{-w/z}}{1+w}$. Comme nous le verrons plus loin, les 3-fonctions apparaissent dans l'étude à l'infini des $E$-fonctions.

---

[6] Prononcer "è" en pensant à Euler à Saint-Petersbourg (suggestion de D. Bertrand).

Nous qualifierons d'*associées* les séries de type $G$, $E$, ou 3 associées à une même suite $\underline{a}$.

D'après 2.1.1, *les séries de type $G$, $E$, 3 ne sont autres que les séries Gevrey de type arithmétique holonomes, d'ordre* $0, -1, +1$, respectivement.

*Exemples de $G$ et $E$-fonctions associées.*

a) Si $\underline{a}$ est la suite constante de valeur 1, on a $f = \frac{1}{1-z}$, $F = e^z$.

b) Si $a_n = n - 1$ pour tout $n$, on a $f = \frac{2z-1}{(1-z)^2}$, $F = (z-1)e^z$.

(Plus généralement, si $f$ est une fonction rationnelle, $F$ est un polynôme exponentiel.)

c) Si $a_0 = 0$, et $a_n = \frac{1}{n}$ pour $n \geq 1$, on a $f = \log(\frac{1}{1-z})$, $F = \frac{e^z - 1}{z} - 1$.

d) Si $a_{2n} = (-1)^n \frac{(\frac{1}{2})_n}{n!}$, et $a_{2n+1} = 0$, on a $f = \frac{1}{\sqrt{1+z^2}}$, $F = J_0(z)$ (fonction de Bessel).

e) Si $a_n = \frac{(\alpha_1)_n \cdots (\alpha_r)_n}{(\beta_1)_n \cdots (\beta_r)_n}$, pour des nombres rationnels $\alpha_i, \beta_i$, avec $\beta_i \notin \mathbb{N}$, (cf. 1.1.2), on a

$$f = {}_{r+1}F_r \begin{pmatrix} \alpha_1, \ldots, \alpha_r, 1 \\ \beta_1, \ldots, \beta_r \end{pmatrix} ; z \end{pmatrix}, F = {}_r F_r \begin{pmatrix} \alpha_1, \ldots, \alpha_r \\ \beta_1, \ldots, \beta_r \end{pmatrix} ; z \end{pmatrix}.$$

En revanche, les séries à coefficients rationnels $\frac{e^z}{1-z}$ et $(1-z)^{\sqrt{2}} + (1-z)^{-\sqrt{2}}$ sont holonomes et Gevrey d'ordre exact 0, mais pas de type arithmétique.

2.3. On peut aussi introduire une version affaiblie de (**G**):

(**G⁻**): Pour tout $\varepsilon > 0$ et pour tout $n$ assez grand,

(**G⁻**)$_{\text{conj}}$: Les conjugués de $a_n$ sont de module inférieur à $(n!)^\varepsilon$,

(**G⁻**)$_{\text{dén}}$: Le dénominateur commun à $a_0, a_1, a_2, \ldots, a_n$ est inférieur à $(n!)^\varepsilon$.

Nous appellerons séries Gevrey de type arithmétique (resp. de type $G$, $E$, ou 3) *au sens large* les séries associées comme ci-dessus à une suite vérifiant (**G⁻**) au lieu de (**G**).

Les $E$-fonctions au sens large sont exactement celles introduites dans l'article original de Siegel [S].[7] Dans les travaux d'approximation diophantienne de l'école de Shidlovskii apparaissent aussi des $E$-fonctions au sens strict de 2.2 (cf. [Sh, p. 406]).

En ce qui concerne les $G$-fonctions, Siegel commence par les définir au sens large, mais précise immédiatement qu'il ne considérera que celles au sens strict (loc. cit. p. 239). On conjecture en fait que les sens stricts et larges coïncident (cf. e.g. [Sh, ibid.]), c'est-à-dire:

CONJECTURE 2.3.1. *Sous* $(\mathbf{H})$, *les conditions* $(\mathbf{G})$ *et* $(\mathbf{G}^-)$ *sont équivalentes.*

Le problème a trait aux dénominateurs, car on a:

PROPOSITION 2.3.2.[8] *Sous* $(\mathbf{H})$, *les conditions* $(\mathbf{G})_{\mathrm{conj}}$ *et* $(\mathbf{G}^-)_{\mathrm{conj}}$ *sont équivalentes.*

*Démonstration.* Soit $K$ comme ci-dessus le corps de nombres engendré par les termes de $\underline{a}$ et par les coefficients d'un polynôme $\phi$ en $z$ et en $\frac{d}{dz}$ non nul qui annule $f = \sum_{n \geq 0} a_n z^n$. Un théorème classique de O. Perron [P] (résultat de base de la théorie des séries Gevrey, voir [R2]) nous place alors dans l'alternative suivante: Ou bien $f$ converge, ou bien il existe un nombre rationnel strictement positif $\eta$ (c'est l'inverse d'une pente du polygône de Newton de $\phi$ en 0) tel que $\underline{a}$ soit sujette à une estimation du type $|a_n| < \kappa(n!)^\eta C^n$, mais à aucune estimation semblable pour un exposant inférieur à $\eta$. On en déduit, en considérant tour à tour chaque plongement complexe $\iota$ de $K$, que $(\mathbf{G}^-)_{\mathrm{conj}}$ entraîne que $f$ définit un germe de fonction analytique en 0 via $\iota$, et ceci équivaut à $(\mathbf{G})_{\mathrm{conj}}$.

*Les séries Gevrey de type arithmétique au sens large n'interviendront plus dans ce volet.*

2.4. On remarque que si la suite $\underline{a}$ vérifie $(\mathbf{G})_{\mathrm{conj}}$, le rayon de convergence de $f = f_{\underline{a}}$ est au moins $C^{-1}$ (relativement à un plongement complexe de $K$ déterminé quelconque) et $F = F_{\underline{a}}$ vérifie $|F(z)| = O(e^{C|z|})$. On peut donc définir la transformée de Laplace $F^+(z) = \int_0^\infty F(w)e^{-wz}dw$ sur le demi-plan Re $(z) > C$ (on pourrait du reste remplacer l'intervalle d'intégration par n'importe quelle demi-droite d'origine 0). Un calcul standard montre alors que

---

[7] L'hypothèse d'holonomie figure dans sa définition (p. 223), mais en a été écartée ultérieurement dans les travaux de l'école russe.

[8] Voir aussi le dernier paragraphe de [R2].

$f(z) = \frac{1}{z}F^+(\frac{1}{z})$, ce qu'on écrira aussi $f(z) = \int_0^\infty F(zw)e^{-w}dw$, ou bien

$$(2.4.1) \qquad\qquad F^+(z) = \frac{1}{z}f\left(\frac{1}{z}\right).$$

Les propriétés suivantes de la transformation de Laplace sont bien connues (voir par exemple [DiP, I.2]):

$$(2.4.2) \qquad\qquad \frac{d}{dz}F^+ = (-zF)^+,$$

$$(2.4.3) \qquad\qquad z \cdot F^+ = \left(\frac{d}{dz}F\right)^+ + F(0),$$

$$(2.4.4) \qquad\qquad \frac{1}{z} \cdot F^+ = \left(\int_0^z F\right)^+,$$

$$(2.4.5) \qquad\qquad F^+(z-a) = (e^{az}F)^+,$$

$$(2.4.6) \qquad\qquad F^+(z/a) = a(F(az))^+.$$

En fait, ces propriétés sont purement formelles: en associant à la suite $\underline{a}$ les séries $f$ et $F \in K[[z]]$ comme ci-dessus, et en définissant la série $F^+ \in \frac{1}{z}K[[\frac{1}{z}]]$ par (2.4.1), on a *ipso facto* les identités (2.4.2) à (2.4.6) (avec $F(0) = a_0$).

On a vu d'autre part, sous la condition (**H**), qu'il existe $\phi \in K[z, \frac{d}{dz}]$ d'ordre $\mu$ annulant $f$. Par changement de variable $z \to \frac{1}{z}$, il existe donc aussi $\Psi \in k[z, \frac{d}{dz}]$ de même ordre annulant $F^+$. Rappelons la transformation de Fourier-Laplace $\triangledown$ des opérateurs différentiels (ou plus généralement des $K[z, \frac{d}{dz}]$-modules) est l'automorphisme d'ordre 4 de $K[z, \frac{d}{dz}]$ défini par

$$z \to -\frac{d}{dz}, \frac{d}{dz} \to z;$$

en notant $\overline{\triangledown}$ la transformation $\triangledown$ suivie de la symétrie par rapport à l'origine, on a donc

$$\overline{\triangledown}\ \triangledown = \triangledown\ \overline{\triangledown} = \mathrm{id}.$$

Posons $\nu = \deg_z \Psi$. La multiplication à gauche de $\Psi$ par $\frac{d^\nu}{dz^\nu}$ neutralise le terme $F(0)$ dans l'application itérée de (2.4.3), et il découle alors de (2.4.2) et (2.4.3) que $\overline{\triangledown}(\frac{d^\nu}{dz^\nu}\Psi)F = 0$. Comme

$$\overline{\triangledown}(\frac{d^\nu}{dz^\nu}\Psi)F = (-z)^\nu\ \overline{\triangledown}(\Psi)F,$$

on voit que $F$ vérifie l'équation différentielle d'ordre $\nu$,

$$(2.4.7) \qquad\qquad \Theta F = 0, \quad \text{avec } \Phi = \overline{\triangledown}(\Psi).$$

Réciproquement, si $F$ vérifie l'équation différentielle $\Phi F = 0$, et si $\mu = \deg_z \Psi$, alors $F^+$ vérifie l'équation différentielle $\frac{d^\mu}{dz^\mu}\triangledown(\Phi)F^+ = 0$, i.e. $\triangledown(\Phi)F^+$ est un polynôme de degré $< \mu$, en général non nul. On retrouve par là en particulier le fait que toute $E$-fonction vérifie une équation différentielle à coefficients polynômiaux. Nous préciserons plus loin les propriétés d'une telle équation.

2.5.  On déduit de 1.1.1 et 2.1.1 que les $G$-fonctions forment une $\mathbb{Q}[z]$-sous-algèbre de $\mathbb{Q}[[z]]$ stable par différentiation; de même pour les $E$-fonctions. Toutefois, on prendra garde que la bijection $\overline{\mathbb{Q}}$-linéaire "série associée"

$$\{\text{G−fonctions}\} \quad \leftrightarrow \quad \{\text{E−fonctions}\}$$
$$f \qquad \leftrightarrow \qquad F$$

ne respecte pas le produit de Cauchy.

De ce que les $E$-fonctions sont entières d'ordre 1, il est facile de déduire que *les unités de l'anneau des $E$-fonctions sont de la forme $be^{az}$*, où $a$ et $b$ sont des nombres algébriques ($b \neq 0$). Il est beaucoup plus délicat de montrer que *les unités de l'anneau des $G$-fonctions sont les fonctions algébriques* holomorphes et non nulles en 0 (cf. [A, p. 124]); nous n'aurons pas à faire usage de ce résultat dans cet article.

Pour toute fonction rationnelle $r(z)$ à coefficients algébriques et s'annulant en 0 et pour toute $G$-fonction $f(z)$, $f(r(z))$ est encore une $G$-fonction.[9]  Il n'en va pas du tout de même pour les $E$-fonctions, pour lesquelles les seuls changements de coordonnées permis sont les homothéties.

Enfin, puisque les conditions (**G**) et (**H**) sont stables par produit terme à terme, tout produit de Hadamard de $G$-fonctions est une $G$-fonction.

*Remarque.* Tout $y \in \text{NGA}\{z\}_s$ annulé par un élément non nul de $\mathbb{C}[z, \frac{d}{dz}]$ est aussi annulé par un élément non nul de $\overline{\mathbb{Q}}[z, \frac{d}{dz}]$. C'est clair si $y \in \overline{\mathbb{Q}}\{z\}_s^A$; en général, on écrit

$$y = \sum_{\alpha, k} \lambda_{\alpha, k, l} z^{\alpha}(\log^k z) y_{\alpha, k, l}, \quad y_{\alpha, k, l} \in \overline{\mathbb{Q}}\{z\}_s^A,$$

les $\alpha$ étant dans des classes distinctes modulo $\mathbb{Z}$, et on voit que les $y_{\alpha, k, l}$ sont holonomes au sens ci-dessus.  On conclut par l'existence[10] de multiples communs dans $\overline{\mathbb{Q}}[z, \frac{d}{dz}]$.

## 3. $G$-opérateurs différentiels (compendium)

3.1.  Soit $K$ un corps de nombres et soit $\phi = \sum_{j=0}^{\mu} Q_j(z)\frac{d^j}{dz^j} \in k[z, \frac{d}{dz}]$ un opérateur différentiel. On lui associe une suite d'opérateurs $(\phi_m)_{m \geq 1}$ de la manière suivante: $\phi_m$ est l'unique élément de $K[z, \frac{d}{dz}]$ divisible à droite par $\phi$

---

[9] En effet cette dernière vérifie clairement une équation différentielle à coefficients polynômiaux, et elle converge à l'origine (relativement à tout plongement complexe du corps des coefficients), ce qui établit les conditions (**H**) et (**G**)$_{\text{conj}}$ pour la suite de ses coefficients de Taylor; on établit (**G**)$_{\text{dén}}$ (resp. (**G**$^-$)$_{\text{dén}}$) en remarquant qu'il existe un entier $N > 1$ tel que le développement de Taylor à l'origine de $r(Nz)$ soit à coefficients entiers algébriques.

[10] Condition de Ore, cf. [Bj, I, 8.4].

et de la forme

$$(3.1.1) \qquad \phi_m = \frac{1}{m!} Q_\mu(z)^m \frac{d^{m+\mu-1}}{dz^{m+\mu-1}} + \sum_{j=0}^{\mu-1} Q_{m,j}(z) \frac{d^j}{dz^j} \,.$$

On dit que $\phi$ est *de type G*, ou est un *G-opérateur* s'il satisfait à la condition suivante introduite par Galochkin [G]:

($\boldsymbol{\mathcal{G}}$) Il existe une constante $\mathcal{C} > 0$ telle que pour tout $n \geq 1$, le dénominateur commun aux coefficients des $Q_{m,j}(z)$, pour $m \leq n$ et $0 \leq j \leq \mu - 1$, est inférieur à $\mathcal{C}^n$.

Cette appellation est justifiée par la remarque (quasi immédiate) que toute solution dans $K[[z-a]]$ d'un $G$-opérateur en un point $a$ tel que $Q_\mu(a) \neq 0$ est une $G$-fonction (en la variable $z - a$).

On constate sans difficulté que si $\phi$ vérifie ($\boldsymbol{\mathcal{G}}$), il en est de même de tout diviseur à droite de $\phi$; il en est aussi de même de l'opérateur $\Psi$ considéré en 2.4, ainsi que de l'opérateur obtenu à partir de $\phi$ en faisant un changement de coordonnées $z \to az$, $a \in K^*$. (En fait, la notion de $G$-opérateur est en un sens évident stable par changement rationnel de coordonnée, et peut se définir sur toute variété algébrique définie sur $\overline{\mathbb{Q}}$; cf. [ABa].)

On *conjecture* que les $G$-opérateurs sont *d'origine géométrique*, c'est-à-dire produits de facteurs d'opérateurs de Picard-Fuchs (contrôlant la variation de cohomologie dans une famille à un paramètre de variétés algébriques définie sur $\overline{\mathbb{Q}}$).

Rappelons les quatre résultats de base concernant ces opérateurs.

THÉORÈME 3.2 (Chudnovsky).     *Soit $f \in K[[z]]$ une G-fonction, et soit $\phi \in K[z, \frac{d}{dz}]$ un opérateur non nul d'ordre minimal tel que $\phi f = 0$. Alors $\phi$ vérifie ($\boldsymbol{\mathcal{G}}$).*

Ce résultat est démontré dans [CC]; voir aussi [A, VI].

Les deux résultats suivants mettent en jeu les *rayons de solubilité générique* $R_v(\phi, 1)$, ainsi définis: pour toute place finie $v$ de K, de caractéristique résiduelle $p = p(v)$, normalisons la valeur absolue v-adique par $|p|_v = p^{-1}$; pour tout $r > 0$, $R_v(\phi, r)$ est le rayon de convergence, limité supérieurement à $r$ par convention, d'une base de solutions de $\phi$ au point générique $v$-adique de valeur absolue $r$.

THÉORÈME 3.3.     *La condition ($\boldsymbol{\mathcal{G}}$) équivaut à $\prod_v R_v(\phi, 1) \neq 0$.*

Ce résultat est démontré dans [A, IV5.2] (la normalisation des valeurs absolues y est différente, mais cela n'a pas d'incidence sur le résultat).

On déduit par exemple de ce critère qu'un produit $\phi \cdot \phi'$ de $G$-opérateurs en est un, de même que l'adjoint $\phi^* = \sum_{j=0}^{\mu}(-1)^j \frac{d^j}{dz^j} \cdot Q_j(z)$; cf. [A, IV4.4].

**THÉORÈME 3.4** (Katz). *Supposons que $R_v(\phi, 1) > p(v)^{-1/(p(v)-1)}$ pour toute place $v$ au-dessus d'un ensemble de nombres premiers de densité 1 (c'est le cas en particulier si $\prod_v R_v(\phi, 1) \neq 0$). Alors $\phi$ est fuchsien, i.e. n'a que des singularités régulières à distance finie et à l'infini. De plus les exposants en chaque singularité sont rationnels.*

Voir [A, IV5.2] ou [DGS, III].

En combinant ces trois résultats (Chudnovsky, André, Katz), on obtient qu'*une équation différentielle d'ordre minimal satisfaite par une $G$-fonction est toujours fuchsienne*, à exposants rationnels. En tout point $\zeta$ (resp. en l'infini), elle admet donc une base de solutions de la forme

$$(3.4.1)_\zeta \qquad \left( f_1^{(\zeta)}(z-\zeta), \ldots, f_\mu^{(\zeta)}(z-\zeta) \right) \cdot (z-\zeta)^{C_\zeta},$$

$$(3.4.2)_\infty \qquad (\text{resp. } \left( f_1^{(\infty)}\left(\frac{1}{z}\right), \ldots, f_\mu^{(\infty)}\left(\frac{1}{z}\right) \right) \cdot \left(\frac{1}{z}\right)^{C_\infty}),$$

où les $f_j^{(\zeta)}, f_j^{(\infty)}$ sont des séries formelles à coefficients dans $K$, et $C_\zeta, C_\infty$ désignent des matrices carrées d'ordre $\mu$ triangulaires supérieures à coefficients dans $\mathbb{Q}$.

**THÉORÈME 3.5** (Sous ($\boldsymbol{\mathcal{G}}$)). *Les séries $f_j^{(\infty)}$ et $f_j^{(\zeta)}$ (pour tout $\zeta \in \overline{\mathbb{Q}}$) sont des $G$-fonctions.*

Ce résultat est démontré dans [A, V]. (C'est non banal seulement si $\zeta$ est une singularité).

3.6. Techniquement, il est plus commode de travailler avec des modules différentiels (sur $\overline{\mathbb{Q}}(z)$) plutôt qu'avec des opérateurs différentiels ([A, IV]). Les définitions sont compatibles au sens où $\overline{\mathbb{Q}}(z)[\frac{d}{dz}]/\overline{\mathbb{Q}}(z)[\frac{d}{dz}]\phi$ est de type $G$ si et seulement si $\phi$ l'est. On dit aussi qu'un $\overline{\mathbb{Q}}[z, \frac{d}{dz}]$-module holonome (i.e. fini et de torsion) est de type $G$ si sa fibre générique l'est.

Le fait que la somme directe et le produit tensoriel de modules différentiels de type $G$ sont de type $G$ ([A, IV4.4]) entraîne le point i) du lemme suivant:

**LEMME 3.6.1.** i) *Soit $y_1$ (resp. $y_2$) une solution d'un $G$-opérateur (resp. d'un autre $G$-opérateur) dans une $\overline{\mathbb{Q}}[z]$-algèbre différentielle convenable. Alors $y_1 + y_2$ est solution d'un $G$-opérateur. De même pour $y_1 \cdot y_2$.*

ii) *Soient $\phi_1$ et $\phi_2$ deux $G$-opérateurs, et soit $\phi$ un multiple commun à gauche de $\phi_1$ et $\phi_2$, d'ordre minimal. Alors $\phi$ est un $G$-opérateur.*

Pour le point i), considérer le sous-module différentiel

$$\overline{\mathbb{Q}}(z)\left[\frac{d}{dz}\right](y_1 + y_2) \quad (\text{resp. } \overline{\mathbb{Q}}(z)\left[\frac{d}{dz}\right]y_1 \cdot y_2)$$

de

$$\overline{\mathbb{Q}}(z)\left[\frac{d}{dz}\right]y_1 \oplus \overline{\mathbb{Q}}(z)\left[\frac{d}{dz}\right]y_2 \quad (\text{resp. } \overline{\mathbb{Q}}(z)\left[\frac{d}{dz}\right]y_1 \otimes_{\overline{\mathbb{Q}}(z)} \overline{\mathbb{Q}}(z)\left[\frac{d}{dz}\right]y_2).$$

Le point ii) découle de i) en utilisant le fait que tout multiple à gauche par un élément de $\overline{\mathbb{Q}}[z]$, et tout diviseur à droite, d'un $G$-opérateur, sont des $G$-opérateurs).

Vu que $z^\alpha \log^k z$ est solution d'un $G$-opérateur (d'ordre $k+1$) lorsque $\alpha$ est rationnel, ce lemme permet d'étendre le théorème de Chudnovsky aux éléments de $\mathrm{NGA}\{z\}_0$:

  *Tout élément holonome de $\mathrm{NGA}\{z\}_0$ est solution d'un $G$-opérateur.*

En combinant ceci au théorème 3.5, on en déduit le *théorème de permanence* énoncé dans l'introduction (et le cas $s = 0$ du théorème de pureté).

## 4. $E$-opérateurs différentiels

4.1. Nous appellerons *$E$-opérateur* le transformé de Fourier-Laplace d'un $G$-opérateur. Il revient au même de dire que $\Phi$ est un $E$-opérateur si et seulement si son transformé de Fourier-Laplace est un $G$-opérateur.

En effet, si $\overline{\Psi}$ désigne l'opérateur déduit de $\Psi := \overline{\mathscr{F}}\Phi$ en appliquant la symétrie par rapport à l'origine, on a $\Phi = \overline{\mathscr{F}}\Psi = \overline{\mathscr{F}\Psi}$; or $\overline{\Psi}$ est un $G$-opérateur si et seulement si $\Psi$ en est un. Un produit $\Phi \cdot \Phi'$ de $E$-opérateurs est un $E$-opérateur, de même que tout diviseur à droite de $\Phi$ et que l'adjoint $\Phi^*$; $\Phi$ et $\Phi'$ admettent un multiple commun à gauche qui est un $E$-opérateur (propriété de Ore à gauche). Ceci résulte des propriétés analogues pour les $G$-opérateurs, et de ce que $(\overline{\mathscr{F}\Psi})^* = \overline{\mathscr{F}}(\Psi^*)$ (cf. [M, V3.6]).

Le nom de "$E$-opérateur" est motivé par le résultat suivant.

THÉORÈME 4.2.   *Toute $E$-fonction $F$ est solution d'un $E$-opérateur $\Phi$; de plus, on peut choisir $\Phi$ de degré en $z$ égal au minimum des degrés en $z$ des opérateurs différentiels annulant $F$. Plus généralement, tout élément holonome de $\mathrm{NGA}\{z\}_{-1}$ est solution d'un $E$-opérateur.*

En effet, soient $f$ la $G$-fonction associée à $F$ et $\phi$ un opérateur différentiel non nul d'ordre minimal annulant $f$. Alors $\phi$ est un $G$-opérateur (3.2), et il en est donc de même de $\Psi$ (cf. 2.4) et de $\overline{\Psi}$. Donc $\Phi := \overline{\mathscr{F}}\Psi = \overline{\mathscr{F}\Psi}$ est un $E$-opérateur non nul, et annule $F$ (cf. 2.4.7).

Pour démontrons la seconde assertion, on peut supposer que cet opérateur $\Phi = \overline{\bigtriangledown}\Psi$ n'est pas de degré minimal en $z$ parmi les opérateurs différentiels annulant $F$. Montrons qu'alors tout opérateur différentiel $\Phi'$ annulant $F$, de degré minimal en $z$ (noté $\mu_0$), est un $E$-opérateur. En effet, $(\overline{\bigtriangledown}\Phi')F$ est un polynôme de degré $< \mu_0$. On en déduit d'une part que $\overline{\bigtriangledown}\Phi'$ est un opérateur d'ordre minimal tel que $\frac{d^{\mu_0}}{dz^{\mu_0}}(\overline{\bigtriangledown}\Phi') \cdot F = 0$, et d'autre part que $\Psi$ divise à droite le produit $\Psi'$ d'un opérateur différentiel d'ordre 1 et de $\overline{\bigtriangledown}\Phi'$. Comme ord $\Psi \neq \overline{\bigtriangledown}\Phi'$ par hypothèse, et puisque $\frac{d^{\mu_0}}{dz^{\mu_0}}\Psi \cdot F = 0$, on en déduit ord $\Psi =$ ord $\Psi'$. Dès lors $\Psi'$, et par suite $\overline{\bigtriangledown}\Phi'$, sont des $G$-opérateurs, donc $\Phi'$ est un $E$-opérateur.

Passons au cas d'un élément holonome $\sum_{\alpha,k}\lambda_{\alpha,k,l}z^{\alpha}(\log^k z)y_{\alpha,k,l}(z)$ de $\mathrm{NGA}\{z\}_{-1}$. En utilisant la propriété de Ore à gauche pour les $E$-opérateurs, on se ramène aisément au cas d'un produit $y(z) = z^{\alpha}(\log^k z)F(z)$, où $F$ est une $E$-fonction (annulée par un $E$-opérateur $\Phi$), et où $\alpha$ est positif; puis, en considérant $\Phi y$ et en raisonnant par récurrence sur $k$, au cas où $k = 0$, i.e. $y = z^{\alpha}F$. On observe alors que la transformée de Laplace $y$ n'est autre que

$$y^+(z) = \Gamma(\alpha+1)z^{-\alpha-1}\left\{\left(1-\frac{1}{z}\right)^{-\alpha} * zF^+(z)\right\},$$

où $*$ désigne le produit de Hadamard de séries en $\frac{1}{z}$. Ainsi $\{(1-\frac{1}{z})^{-\alpha} * zF^+(z)\}$ est une $G$-fonction en $\frac{1}{z}$, et $y^+(z)$ est solution d'un $G$-opérateur (cf. 2.5, 3.2).

*Remarque.* On déduit de 4.2. qu'*un opérateur d'ordre minimal $\Theta$ annulant $F$ admet un multiple à gauche de la forme $Q(z)\Phi$, où $Q$ est un polynôme et $\Phi$ un $E$-opérateur. On prendra toutefois garde qu'en général ni $Q(z)\Phi$ ni $\Theta$ ne sont eux-mêmes des $E$-opérateurs.* C'est le cas de l'exemple 2.3.b ci-dessus ($F = (z-1)e^z$), pour lequel on peut prendre:

$$\phi = (z-1)(2z-1)\frac{d}{dz} + 3z - 1, \quad \Psi = -(z-1)(z-2)\frac{d}{dz} + z - 3,$$

$$\Phi = \overline{\bigtriangledown}\Psi = z\frac{d^2}{dz^2} + (1-3z)\frac{d}{dz} + 2z, \quad \Theta = (z-1)\frac{d}{dz} - z.$$

On a $(z-1)\Phi = (z\frac{d}{dz} - 2z + 1) \cdot \Theta$, donc $Q(z) = z - 1$ convient, mais $\overline{\bigtriangledown}\Theta = -(1+z)\frac{d}{dz} + z - 1$ n'est pas un $G$-opérateur. Notons d'ailleurs que $\Phi$ n'est pas de degré minimal en $z$: on a

$$\left(\frac{d}{dz} - 1\right) \cdot \theta = (z-1)\left(\frac{d^2}{dz^2} - 2 + 1\right),$$

et $\Phi' := \frac{d^2}{dz^2} - 2\frac{d}{dz} + 1$ est un autre $E$-opérateur d'ordre minimal annulant $F$, indépendant de $\Phi$.

THÉORÈME 4.3 (structure des $E$-opérateurs). *Soit $\Phi \in K[z, \frac{d}{dz}]$ un $E$-opérateur, d'ordre noté $\nu$. Alors:*

i) $\Phi$ *n'a que deux singularités*: $0$ *et* $\infty$.

ii) *L'origine est une singularité régulière. Les exposants de $\Phi$ en $0$ sont rationnels; exceptés peut-être ceux entiers, ce sont, modulo $\mathbb{Z}$ (et comptés sans multiplicité), les exposants à l'infini du $G$-opérateur $\overline{\ominus}\Phi$.*

iii) *Il existe une base de solutions en $0$ de l'équation $\Phi F = 0$ de la forme*

$$(F_1(z), \ldots, F_\nu(z)) \cdot z^{\Gamma_0}$$

*où les $F_j$ sont des $E$-fonctions à coefficients dans $K$, et $\Gamma_0$ désigne une matrice carrée d'ordre $\nu$ triangulaire supérieure à coefficients dans $\mathbb{Q}$.*

iv) *L'infini est une singularité en général irrégulière, de pentes $\in \{0, 1\}$. L'ensemble des exposants de Turrittin non entiers de $\Phi$ en $\infty$ coïncide, modulo $\mathbb{Z}$, avec l'ensemble des exposants non entiers de $\overline{\ominus}\Phi$ aux singularités à distance finie.*

v) *Il existe une base de solutions en $\infty$ de l'équation $\Phi F = 0$ de la forme*

$$\left(\mathfrak{f}_1\left(\frac{1}{z}\right), \ldots, \mathfrak{f}_\nu\left(\frac{1}{z}\right)\right) \cdot \left(\frac{1}{z}\right)^{\Gamma_\infty} \cdot e^{-\Delta z}$$

*où les $\mathfrak{f}_j$ sont des $3$-fonctions à coefficients dans l'extension finie de $K$ obtenue par adjonction des singularités à distance finie de $\overline{\ominus}\Phi$, où $\Delta$ est la matrice diagonale ayant pour coefficients diagonaux les singularités à distance finie de $\overline{\ominus}\Phi$ (comptées avec multiplicité), et où $\Gamma_\infty$ désigne une matrice carrée triangulaire supérieure à coefficients dans $\mathbb{Q}$ (ayant pour éléments diagonaux les exposants de Turrittin), qui commute à $\Delta$.*

La démonstration occupera les sections 5.1 à 5.6. Nous priviligierons une approche formelle, qui permettrait de travailler sur $\overline{\mathbb{Q}}(\Gamma^{(k)}(a))_{k\in\mathbb{N}, a\in\mathbb{Q}}$ plutôt que sur $\mathbb{C}$ si l'on voulait (en fait les $\Gamma^{(k)}(a)$ n'interviennent que via leurs "relations de distribution"). Toutefois nous préciserons aussi par voie analytique les liens entre $\Phi$ et $\overline{\ominus}\Phi$.

*Exemples.* a) L'opérateur à coefficients constants $\Phi = \frac{d^2}{dz^2} - 2\frac{d}{dz} + 1$ considéré plus haut admet les $E$-fonctions $e^z$ et $ze^z$ comme base de solutions; $\overline{\ominus}\Phi$ est l'opérateur d'ordre nul $(z+1)^2$, qui admet $-1$ comme singularité triviale double.

b) Par le point i) du théorème, on voit que les seuls opérateurs différentiels irréductibles qui sont à la fois de type $G$ et $E$ sont, à multiplication près par un monôme en $z$, ceux de la forme $z\frac{d}{dz} - \alpha$ avec $\alpha \in \mathbb{Q}^*$, et $\frac{d}{dz}$.

c) Equation différentielle de Whittaker [WW], [Bu]. Soient $k$ et $m$ deux nombres rationnels, et considérons *l'opérateur de Whittaker*

$$\Phi_{k,m} = z^2\frac{d^2}{dz^2} + \left(-\frac{z^2}{4} + kz + \frac{1}{4} - m^2\right).$$

On calcule

$$\overline{\overline{\mathcal{7}}}\,\Phi_{k,m} = \left(z - \frac{1}{2}\right)\left(z + \frac{1}{2}\right)\frac{d^2}{dz^2} + (4z + k)\frac{d}{dz} + \left(\frac{9}{4} - m^2\right),$$

qui est un $G$-opérateur, car c'est un opérateur de type hypergéométrique à exposants rationnels (correspondant au schéma de Riemann

$$\mathfrak{P}\begin{pmatrix} -\frac{1}{2} & +\frac{1}{2} & \infty \\ 0 & 0 & -m + \frac{3}{2} \quad ; z \\ -k - 1 & +k - 1 & +m + \frac{3}{2} \end{pmatrix}).$$

On rappelle que si $m$ n'est pas demi-entier, une base de solutions de $\Phi_{k,m}$ en 0 est donnée par

$$M_{k,\pm m}(z) = z^{\frac{1}{2}\pm m}e^{-\frac{1}{2}z}{}_1F_1\left(\frac{1}{2}\pm m - k; \pm 2m + 1; z\right);$$

une autre base est donnée par les fonctions de Whittaker

$$W_{k,m}(z) = \frac{\Gamma(-2m)}{\Gamma\left(\frac{1}{2} - m - k\right)}M_{k,m}(z)$$
$$+ \frac{\Gamma(2m)}{\Gamma\left(\frac{1}{2} + m - k\right)}M_{k,-m}(z) \quad \text{et} \quad W_{-k,m}(-z).$$

Le développement asymptotique de $W_{k,m}(z)$ à l'infini est

$$W_{k,m}(z) \approx z^k e^{-\frac{1}{2}z}{}_2F_0\left(\frac{1}{2} - m - k, \frac{1}{2} + m - k; \left(-\frac{1}{z}\right)\right)$$

pour $|\arg(z)| \leq \pi - \varepsilon$. On constate que ${}_1F_1(\frac{1}{2} \pm m - k; \pm 2m + 1; z)$ est une $E$-fonction tandis que ${}_2F_0(\frac{1}{2} - m - k, \frac{1}{2} + m - k; z)$ est une 3-fonction.

d) $E$-opérateur associé à l'équation différentielle de Lamé [WW], [D]. Soient $n$ un nombre rationnel, $B, e_1, e_2, e_3$, quatre nombres algébriques ($e_i \neq e_j$), et considérons *l'opérateur de Lamé*

$$\Lambda_{n,B} = 4(z - e_1)(z - e_2)(z - e_3)\left\{\frac{d^2}{dz^2} + \frac{1}{2}\left(\frac{1}{z - e_1} + \frac{1}{z - e_2} + \frac{1}{z - e_3}\right)\frac{d}{dz}\right.$$
$$\left. - (n(n+1)z + B)\right\}.$$

C'est un opérateur fuchsien à exposants rationnels (correspondant au schéma de Riemann

$$\mathfrak{P}\begin{pmatrix} e_i & \infty \\ 0 & -\frac{n}{2}; z \\ \frac{1}{2} & \frac{n+1}{2} \end{pmatrix}).$$

Si $B = B_n^m$ est l'une des $2n+1$ valeurs classiques pour lesquelles $\Lambda_{n,B}$ admet une fonction de Lamé solution, il est connu que $\Lambda_{n,B}$ est un $G$-opérateur. L'exemple du $E$-opérateur du troisième ordre $\overline{\triangledown}(\Lambda_{n,B_n^m})$ montre qu'il est nécessaire de se limiter aux exposants de Turrittin *non entiers* au point iv) de 4.3.

Notons par ailleurs que si n est un entier pair négatif, $\overline{\triangledown}(\Lambda_{n,B_n^m})$ admet une $E$-fonction solution (transformée de Laplace inverse d'une fonction de Lamé); il me paraît douteux qu'elle soit hypergéométrique au sens de Siegel en général (cf. [Sh, p. 184]).

Compte tenu de la remarque suivant 4.2, le théorème 4.3 entraîne:

COROLLAIRE 4.4.    *Soit $F$ une $E$-fonction, et soit $\Theta \in K[z, \frac{d}{dz}]$ un opérateur non nul d'ordre minimal annulant $F$. Alors $\Theta$ n'a que deux singularités non triviales*: $0$ *et* $\infty$.

*L'origine est une singularité régulière, à exposants rationnels; plus précisément, il existe une base de solutions en 0 de l'équation $\Theta F = 0$ de la forme $(F_1(z), \ldots, F_\lambda(z)) \cdot z^\Gamma$, où les $F_j$ sont des $E$-fonctions à coefficients dans $K$ (et $F_1 = F$), et où $\Gamma$ désigne une matrice carrée triangulaire supérieure à coefficients dans $\mathbb{Q}$.*

*Remarque.* Le fait que $\Theta$ n'ait que deux singularités non triviales entraîne que l'algèbre de Lie de son groupe de Galois différentiel global se calcule arithmétiquement selon la conjecture de Katz [K], [B]: c'est la plus petite sousalgèbre de Lie de l'algèbre de Lie des matrices carrées d'ordre celui de $\Theta$, dont la réduction modulo presque tout nombre premier $p$ contient la $p$-courbure de $\Theta$.

4.5. *Application à un problème de Shidlovskii.* Dans [Sh, p. 184], Shidlovskii remarque que toutes les $E$-fonctions solutions d'une équation différentielle homogène (linéaire) d'ordre 1 sont de la forme $Q(z)e^{\zeta z}$, où $Q$ est un polynôme, et signale que la détermination des $E$-fonctions solutions d'une équation différentielle *inhomogène d'ordre 1* est une question ouverte. Le corollaire cidessus permet d'y apporter une réponse, en considérant l'équation homogène associée.

On trouve que ce sont les fonctions entières de la forme

$$F = q(z)z^{-\alpha}e^{\zeta z}\int^z r(t)t^{\alpha-1}e^{-\zeta t}$$

où $\alpha$ est rationnel, $\zeta$ est algébrique, et où $q(z) \in \overline{\mathbb{Q}}[z]$, $r(z) \in \overline{\mathbb{Q}}(z)$ (soumis aux conditions adéquates pour que $F$ soit entière). Le prototype en est

$$q(z)_1F_1(1;\alpha;z) + p(z), \text{ avec } p(z), q(z) \in \overline{\mathbb{Q}}[z];$$

il resterait à déterminer si on les obtient toutes ainsi.

Le lien tissé entre 3-fonctions et $E$-fonctions par le théorème 4.3 est renforcé par le résultat suivant, pendant de 4.2:

THÉORÈME 4.6.     *Soit $\mathfrak{f}(z)$ une 3-fonction (ou plus généralement un élément holonome de $\mathrm{NGA}\{z\}_{+1}$). Alors $\mathfrak{f}(\frac{1}{z})$ est solution d'un $E$-opérateur.*

Par exemple $\sum_{n\geq0}(-1)^n n!(\frac{1}{z})^n$ est solution du $E$-opérateur $\Phi = z\frac{d^2}{dz^2} + (1-z)\frac{d}{dz} - 1$ (de $G$-opérateur associé $\overline{\forall}\Phi = z((1+z)\frac{d}{dz}+1)$).

Je conjecture la caractérisation $p$-adique suivante des $E$-opérateurs (pendant de 3.3):

CONJECTURE 4.7. *Supposons que $\Phi \in K[z, \frac{d}{dz}]$ n'ait de singularité qu'en $0$ et $\infty$. Alors $\Phi$ est un $E$-opérateur si et seulement si*

$$\prod_v \left(R_v\left(\Phi, p(v)^{-1/p(v)-1}\right)\cdot p(v)^{1/(p(v)-1)}\right) \neq 0.$$

*Remarque* 4.8.     En interprétant la transformée de Fourier-Laplace comme un foncteur $p_{2*}(p_1^*? \otimes e^{-zw})[1]$ en théorie des $D$-modules (cf. [M, app. 2]), les méthodes de [AB] ramènent essentiellement la conjecture à un problème de convergence $p$-adique des séries formelles intervenant dans la décomposition de Turrittin-Levelt à l'origine et à l'infini, pour $p$ assez grand.

4.9. Plus généralement, on définit les $\overline{\mathbb{Q}}[z, \frac{d}{dz}]$-modules de type $E$ comme étant les transformés de Fourier-Laplace des $\overline{\mathbb{Q}}[z, \frac{d}{dz}$-modules holonomes de type $G$ (cf. 3.6). Les définitions sont compatibles au sens où $\overline{\mathbb{Q}}[z, \frac{d}{dz}]/\overline{\mathbb{Q}}[z, \frac{d}{dz}]\Phi$ est de type $E$ si et seulement si $\Phi$ l'est. On prendra garde que cette propriété ne dépend pas uniquement de la fibre générique de $\overline{\mathbb{Q}}[z, \frac{d}{dz}]/\overline{\mathbb{Q}}[z, \frac{d}{dz}]\Phi$, contrairement à ce qui se passe pour le "type $G$".

## 5. Autour de Laplace. Preuve de 4.3 et 4.6

5.1. *Pentes.* Soit $\varphi = \sum_{i=0}^{\mu}\sum_{j=0}^{\nu}a_{i,j}z^j\frac{d^i}{dz^i} \in K[z, \frac{d}{dz}]$; son transformé de Fourier-Laplace est

$$\overline{\forall}\varphi = \sum_{j=0}^{\nu}\sum_{i=0}^{\mu}(-1)^ja_{i,j}\frac{d^j}{dz^j}z^i = \sum_{j=0}^{\nu}\sum_{i=0}^{\mu}b_{j,i}z^i\frac{d^j}{dz^j}.$$

Tous les opérateurs différentiels que nous considérons sont de type exponentiel au sens de [M, XII], c'est-à-dire vérifient $a_{\mu,\nu} \neq 0$ (ce qui équivaut à $b_{\nu,\mu} \neq 0$). En effet nous ne considérons que des opérateurs différentiels réguliers en l'infini, ou facteurs de transformés de Fourier de tels.

Notons d'autre part que si $\varphi$ est régulier en l'infini, on a alors $a_{i\nu} = 0$ pour $i < \mu$, donc $b_{\nu,i} = 0$ pour $i < \mu$, ce qui signifie que $\overline{\forall}\varphi$ n'a de singularités qu'en $0$ et $\infty$. Ceci établit i).

Pour aller plus loin, rappelons que le *polygône de Newton-Ramis $N(\varphi)$* est l'enveloppe convexe dans le plan des demi-droites $\{u \leq i, v = j - i/a_{i,j} \neq 0\}$ [R1], [M, V1]. Toute pente non nulle n'apparaît que sur un côté; on distingue les côtés horizontaux en donnant la pente $-0$ (resp. $+0$) au côté supérieur (resp. inférieur). De la sorte, la partie de pente $> 0$ (resp. $< 0$) est à translation verticale près le polygône de Newton usuel de $\varphi$ en 0 (resp. le symétrisé par rapport à l'axe des $u$ du polygône de Newton usuel de $\varphi$ en $\infty$).

Il est immédiat qu'on passe de $N(\varphi)$ à $N(\overline{\vee}\varphi)$ par la transformation $(u, v) \to (u + v, -v)$; en particulier les pentes $t$ de $N(\varphi)$ devient la pente $\frac{-t}{t+1}$ de $N(\overline{\vee}\varphi)$. On en déduit tout de suite le critère suivant (avec $\Phi = \overline{\vee}\varphi$):

$(\Phi$ *n'a que deux singularités,* 0 *qui est régulière,*

*et l'infini qui est irrégulière de pentes* $\in \{0, 1\})$

$\Leftrightarrow ( \overline{\vee}\Phi$ *est régulier en* 0 *et en l'infini).*

Compte tenu de ce que tout $G$-opérateur est fuchsien, ceci établit les assertions de 4.3 concernant les pentes de $\Phi$.

On peut dire plus: si $\varphi$ est fuchsien (i.e. n'a que des singularités régulières), alors on peut appliquer ci-dessus aux translatés $\varphi_a$ de $\varphi$ (par $z \to z - a$); comme $\overline{\vee}\varphi_a$ n'est autre que le "tordu" de $\overline{\vee}\varphi$ par $e^{az}$, on en conclut que les pentes d'un tel tordu sont toujours dans $\{0, 1\}$, ce qui montre que les facteurs déterminants de $\overline{\vee}\varphi$ (*a priori* constants ou bien polynômiaux de degré $\nu!$ en $z^{1/\nu!}$) sont simplement de la forme $\delta z$. En particulier, d'après Turrittin-Levelt, il existe une base de solutions en $\infty$ de l'équation $\Phi F = 0$ de la forme

$$\left( \hat{f}_1\left(\frac{1}{z}\right), \cdots, \hat{f}_\nu\left(\frac{1}{z}\right) \right) \cdot \left(\frac{1}{z}\right)^{\Gamma_\infty} \cdot e^{-\Delta z}$$

où $\Gamma_\infty$ désigne une matrice sous forme de Jordan (ayant pour éléments diagonaux les exposants de Turrittin), où $\Delta$ est une matrice diagonale commutant à $\Gamma_\infty$ ($\Delta z$ ayant pour éléments diagonaux les facteurs déterminants), et où les $\hat{f}_j$ sont des séries de Laurent à coefficients dans l'extension de $K$ obtenue par adjonction des coefficients de $\Gamma_\infty$ et $\Delta$.

5.2. *Monodromie de et monodromie de* $\overline{\vee}\Phi$. Fixons un plongement complexe de $K$. Notons $\{\zeta_1, \ldots, \zeta_r\}$ l'ensemble des singularités finies de $\varphi$, qu'on suppose fuchsiennes (de même que l'infini). Alors le module différentiel $\mathbb{C}[z, \frac{d}{dz}]/\mathbb{C}[z, \frac{d}{dz}]\varphi$ est déterminé à isomorphisme près par le système local des solutions de $\varphi$ sur $\mathbb{C} - \{\zeta_1, \ldots, \zeta_r\}$.

Pour $\Phi = \overline{\vee}\varphi$, en revanche, le module différentiel $\mathbb{C}[z, \frac{d}{dz}]/\mathbb{C}[z, \frac{d}{dz}]\Phi$ est déterminé à isomorphisme près par le système local des solutions de $\Phi$ sur $\mathbb{C}^*$ *et* par la structure de Stokes à l'infini (cf. [M, IV]). La détermination de ces données est de manière générale l'objet de la théorie de la *transformation de Fourier géométrique*, mais dans le cas exponentiel élémentaire que nous considérons ici, elle est fournie par la recette explicite suivante ([M, XII]).

Pour l'exposer, nous introduirons comme dans loc. cit. une autre variable, soit $w$, et nous mettrons en dualité réelle le plan complexe des $z$ (avec l'orientation usuelle) et le plan complexe des $w$ (avec l'orientation inhabituelle). Nous considérerons $\Phi$ comme un élément de $\mathbb{C}[w, \frac{d}{dw}]$; sa monodromie en 0 tourne donc par convention dans le sens inhabituel.

Soit $\widetilde{\mathcal{O}_z}$ le $\mathbb{C}[z, \frac{d}{dz}]$-module des fonctions holomorphes sur le revêtement universel du plan des $z$ privé d'un disque fermé assez grand, et soit $\widetilde{\mathcal{C}_z}$ son quotient par les fonctions entières d'ordre $\leq 1$. On a les applications standard $\widetilde{\mathcal{O}_z} \overset{\mathrm{Can}}{\underset{\mathrm{Var}}{\rightleftarrows}} \widetilde{\mathcal{C}_z}$, et $T = \mathrm{id} + \mathrm{Var}\,\mathrm{Can}$ est la monodromie autour du disque manquant. De même pour la variable $w$, avec l'orientation inverse.

Ces applications passent aux noyaux de $\varphi$ et de $\Phi$ respectivement:

$$\mathrm{Ker}(\varphi, \widetilde{\mathcal{O}_z}) \overset{\mathrm{Can}_\varphi}{\underset{\mathrm{Var}_\varphi}{\rightleftarrows}} \mathrm{Ker}(\varphi, \widetilde{\mathcal{C}_z}), \quad \mathrm{Ker}(\Phi, \widetilde{\mathcal{O}_w}) \overset{\mathrm{Var}_\Phi}{\underset{\mathrm{Can}_\Phi}{\rightleftarrows}} \mathrm{Ker}(\Phi, \widetilde{\mathcal{O}_w}).$$

Le résultat principal de loc. cit. est d'une part que *la structure de Stokes de $\Phi$ est indexée par* $\{-\zeta_1, \ldots, -\zeta_r\}$ (ce qui se traduit ici par le fait que dans la décomposition de Turrittin-Levelt au voisinage de $w = \infty$, les facteurs déterminants sont les $-\zeta_i w$),[11] et d'autre part que $\mathrm{Var}_\Phi$ *et* $\mathrm{Can}_\Phi$ *s'identifient canoniquement à* $\mathrm{Can}_\varphi$ *et* $\mathrm{Var}_\varphi$ respectivement.

L'isomorphisme canonique $\mathrm{Ker}(\Phi, \widetilde{\mathcal{O}_w}) = \mathrm{Ker}(\varphi, \widetilde{\mathcal{C}_z})$ est donné de la manière suivante. On a une décomposition $\mathrm{Ker}(\varphi, \widetilde{\mathcal{C}_z}) = \oplus_\zeta \mathrm{Ker}(\varphi, \widetilde{\mathcal{C}_{z=\zeta}})$ en espaces de microsolutions de $\varphi$ en chaque singularité à distance finie. On choisit une direction de demi-droite $\theta$ telle que les demi-droites $\delta_{\zeta, \theta}$, de direction $\theta$ issues de $\zeta$ ne se chevauchent pas. Pour $f_\zeta \in \mathrm{Ker}(\varphi, \widetilde{\mathcal{C}_{z=\zeta}})$, on choisit $\tilde{f}_\zeta$ tel que $\mathrm{can}_\zeta \tilde{f}_\zeta = f_\zeta$. Alors la formule

$$\int_{|z - \zeta| = \varepsilon} \tilde{f}(z) e^{-zw} dz + \int_{\delta\zeta, \theta, |z - \zeta| \geq \varepsilon} (\mathrm{var}_\zeta f)(z) e^{-zw} dz$$

définit une solution de $\Phi$ (qui ne dépend pas du choix de $\tilde{f}$) sur un demi-plan $\mathrm{Re}\,(we^{i\theta}) \gg 0$.[12] On obtient ainsi un sous-espace $\mathrm{Ker}(\Phi, \widetilde{\mathcal{O}_w})_\zeta$ de $\mathrm{Ker}(\Phi, \widetilde{\mathcal{O}_w})$, et on a $\mathrm{Ker}(\Phi, \widetilde{\mathcal{O}_w}) = \oplus_\zeta \mathrm{Ker}(\Phi, \widetilde{\mathcal{O}_w})_\zeta$ (décomposition de Stokes).

---

[11] Pour se convaincre qu'on a pris les bons signes, on peut considérer le cas de $\Phi = \frac{d}{dz} + \zeta$.

[12] La procédure est réversible et permet de construire les solutions de $\varphi$ à l'infini à partir des microsolutions de $\Phi$ à l'origine.

Les flèches $\mathrm{can}_\zeta, \mathrm{var}_\zeta$ et $T_\zeta = \mathrm{id} + \mathrm{var}_\zeta \mathrm{can}_\zeta$ associées aux microsolutions de $\varphi$ en $\zeta$ sont reliées à $\mathrm{Can}_\varphi$ et $\mathrm{Var}_\varphi$ par les formules

$$\mathrm{Can}_\varphi = \begin{pmatrix} \mathrm{can}_{\zeta_1} \\ \vdots \\ \mathrm{can}_{\zeta_i} \cdot T_{\zeta_{i-1}} \cdots T_{\zeta_1} \\ \vdots \\ \mathrm{can}_{\zeta_r} \cdot T_{\zeta_{r-1}} \cdots T_{\zeta_1} \end{pmatrix}, \quad \mathrm{Var}_\varphi = (\mathrm{var}_{\zeta_1}, \ldots, \mathrm{var}_{\zeta_r})$$

(dans loc. cit. ce calcul est fait en supposant, pour fixer les idées, que $\mathrm{Im}\,\zeta_1 > \cdots > \mathrm{Im}\,\zeta_r$, de sorte qu'on puisse prendre $\theta = 0$; mais c'est sans importance, car on se ramène immédiatement à ce cas en effectuant une homothétie sur $z$ et l'homothétie inverse sur $w$).

La monodromie $T_{\Phi,0}$ de $\Phi$ tournant dans le sens usuel autour de 0 s'identifie à $(\mathrm{id} + \mathrm{Var}_\Phi \mathrm{Can}_\Phi)^{-1}$, et donc à $(\mathrm{id} + \mathrm{Can}_\varphi \mathrm{Var}_\varphi)^{-1}$, tandis que la monodromie $T_{\varphi,\infty}$, de $\varphi = \overline{\overline{\mathscr{T}}} \Phi$ autour de l'infini s'identifie à $(\mathrm{id} + \mathrm{Var}_\varphi \mathrm{Can}_\varphi)^{-1}$. Comme les valeurs propres non nulles de $\mathrm{Can}_\varphi \mathrm{Var}_\varphi$ coïncident avec celles de $\mathrm{Var}_\varphi \mathrm{Can}\varphi$, on voit en particulier que les valeurs propres de $T_{\Phi,0}$ et de $T_{\varphi,\infty}$, distinctes de 1 coïncident (compte non tenu des multiplicités), ce qui démontre, par voie analytique, l'assertion 4.3.ii).

Avec ces résultats, joints au théorème de Turritin-Levelt, il ne reste plus qu'à établir que les séries $F_i$ de 4.3. iii) sont des $E$-fonctions, tandis que les séries $\mathfrak{f}_i$ de 4.3. v) sont des 3-fonctions, et à caractériser les exposants de Turritin.

5.3. *Calcul opérationnel* (*formulaire*). Rappelons que pour tout $\alpha$ de partie réelle $> -1$, la transformée de Laplace de $z^\alpha$ est

$$(5.3.1) \qquad (z^\alpha)^+ = \Gamma(\alpha+1) z^{-\alpha-1};$$

plus généralement, pour tout entier naturel $k$,

$$(5.3.2) \quad (z^\alpha \log^k z)^+ = \Gamma(\alpha+1) z^{-\alpha-1} \Big( (-1)^k \log^k z + \text{un polynôme}$$
$$\text{de degré } k-1 \text{ en } \log z \Big)$$

Ce polynôme fait intervenir les valeurs en $\alpha$ de la fonction $\Gamma$ et de ses dérivées jusqu'à l'ordre $k$, mais son expression importe peu ici.

Il nous faudra étendre la transformation $^+$ dans le cas d'un nombre complexe $\alpha$ quelconque, de manière à disposer encore des formules du type de celles de 2.4. Ce problème est résolu par le *calcul opérationnel* de la manière suivante (cf. [DiP, II5.3]).

Considérons une somme finie

$$h = \sum_{\alpha,k} h_{\alpha,k} z^\alpha \log^k z,$$

et posons

$$(5.3.3) \qquad \vartheta(z, \alpha, k) = \frac{z^{\alpha+1}}{\alpha+1} \sum_{l=0}^{k} \frac{k!}{(1-k)!} (-\alpha-1)^{-1} \log^{k-1} z \quad \text{si } \alpha \neq -1,$$

$$(5.3.4) \quad \vartheta(z, -1, k) = \frac{\log^{k+1} z}{k+1} \ .$$

Rappelons que la *partie finie* de l'intégrale $\int_0^z h(t)dt$ est la fonction de $z$

$$\text{p.f.} \int_0^z h(t)dt = \lim_{\varepsilon \to 0^+} \left( \sum_{\alpha, k} h_{\alpha, k} \vartheta(\varepsilon, \alpha, k) + \int_\varepsilon^z h(t)dt \right)$$

$$= \sum_{\alpha, k} h_{\alpha, k} \vartheta(z, \alpha, k).$$

Il suit de cette définition que

$$(5.3.5) \ \text{p.f.} \int_0^z \frac{(z-t)^n}{n!} t^\alpha \log^k t \, dt \ = \ \sum_{m=0, m \neq \alpha-1}^{n} \frac{(-1)^m}{m!(n-m)!} \frac{z^{\alpha+n+1}}{m+\alpha+1}$$

$$\times \left( \log^k z + \cdots + \frac{(-1)^k k!}{(m+\alpha+1)^k} \right)$$

plus un terme

$$\frac{(-1)^{\alpha+1}}{(-\alpha-1)!(n+\alpha+1)!} z^{\alpha+n+1} \frac{\log^{k+1} z}{k+1}$$

si $\alpha - 1$ est un entier compris entre $0$ et $n$.

On a pour ces parties finies un calcul intégro-différentiel analogue au calcul usuel, la formule habituelle $g(z) - \int_0^z (\frac{d}{dt} g(t))dt = g(0)$ étant remplacée par

$$h(z) - \text{p.f.} \int_0^z \left( \frac{d}{dt}(h(t))dt \right) = h_{0,0}$$

(constante qu'on note aussi p.f.$(h)(0)$). On peut alors étendre le formalisme de la transformation de Laplace en posant

$$(5.3.6) \qquad h^+ = z^{n+1} \left( \text{p.f.} \int_0^z \frac{(z-t)^n}{n!} h(t)dt \right)^+$$

pour tout $n \geq -\text{Re}\,\alpha - 1$ (cette expression ne dépend pas d'un tel $n$). Par exemple, la formule (5.3.1) vaut pour tout $\alpha$ non entier négatif, tandis que

$$\left( \frac{1}{z} \right)^+ = z \cdot (\log z)^+ = \Gamma'(1) - \log z.$$

La formule (2.4.2) est encore valide, i.e.

$$(5.3.7) \qquad\qquad \frac{d}{dz} h^+ = (-zh)^+,$$

tandis que (2.4.3) est à modifier comme suit

$$(5.3.8) \qquad z \cdot h^+ = \left(\frac{d}{dz}h\right)^+ + \sum_{m \geq 0} (-1)^m \frac{z^m}{m!} \, \text{p.f.}(z^m h)(0)$$

(la somme s'arrête à $m = n$), ce qui redonne (2.4.3) si $h$ n'a pas de terme qui soit une puissance entière strictement négative de $z$.

Compte tenu de (5.3.5) et de la formule

$$(5.3.9) \qquad \sum_{m=0}^{n} \frac{(-1)^m}{m!(n-m)!} \frac{1}{m+\alpha+1} = ((\alpha+1)_{n+1})^{-1}$$

si $-\alpha - 1$ n'est pas un entier compris entre 0 et $n$, on obtient la généralisation suivante de (5.3.2):

$$(5.3.10) \quad (z^\alpha \log^k z)^+ = \Gamma(\alpha+1) z^{-\alpha-1} ((-1)^k \log^k z$$

$$+ \text{ un polynôme de degré } k-1 \text{ en } \log z)$$

$$\text{si } \alpha \text{ n'est pas un entier strictement négatif,}$$

$$= \frac{z^{-\alpha-1}}{(-\alpha-1)!} ((-1)^\alpha \frac{\log^{k+1} z}{k+1}$$

$$+ \text{ un polynôme de degré } k \text{ en } \log z) \text{ sinon.}$$

5.4. *Laplace et Nilsson-Gevrey.* Par complétion formelle $z$-adique (resp. $\frac{1}{z}$-adique), la transformation $^+$ de (5.3.6) se prolonge en deux applications :

$$z^\alpha \mathbb{C}[[z]][\log z] \underset{\leftarrow}{\overset{\rightarrow}{}} \frac{1}{z^{\alpha+1}} \mathbb{C}\left[\left[\frac{1}{z}\right]\right] \left[\log \frac{1}{z}\right]$$

qui vérifient (5.3.7) et (5.3.8). En filtrant par le degré du logarithme, on voit immédiatement que ces applications sont injectives, et même bijectives si $\alpha$ n'est pas entier. Pour $\alpha$ rationnel, ces applications s'étendent respectivement par linéarité en des applications injectives encore notées $^+$, qui vérifient (5.3.7) et (5.3.8):

$$\overline{\mathbb{C}[[z]]}[\log z] \underset{\leftarrow}{\overset{\rightarrow}{}} \overline{\mathbb{C}\left[\left[\frac{1}{z}\right]\right]} \left[\log \frac{1}{z}\right],$$

où $\overline{\mathbb{C}[[z]]}$ et $\overline{\mathbb{C}[[\frac{1}{z}]]}$ désignent les anneaux de séries de Puiseux en $z$ et en $\frac{1}{z}$ respectivement.

PROPOSITION 5.4.1. *Ces applications induisent des applications $\mathbb{C}$-linéaires injectives*

$$\text{NGA}\{z\}_{-1} \underset{\leftarrow}{\overset{\rightarrow}{}} \text{NGA}\{\tfrac{1}{z}\}_0, \text{NGA}\{z\}_0 \underset{\leftarrow}{\overset{\rightarrow}{}} \text{NGA}\left\{\frac{1}{z}\right\}_{+1},$$

*la source de chacune de ces applications coïncidant avec sa préimage dans* $\overline{\mathbb{C}}[[z]][\log z]$ *ou* $\overline{\mathbb{C}}[[\frac{1}{z}]][\log \frac{1}{z}]$ *respectivement.*

*Preuve.* Nous traiterons la question de $\mathrm{NGA}\{z\}_{-1}\mathrm{NGA}\{\frac{1}{z}\}_0$ et de $\mathrm{NGA}\{\frac{1}{z}\}_1 \to \mathrm{NGA}\{z\}_0$, les deux cas restants étant similaires.

Par $\overline{\mathbb{Q}}$-linéarité, on se ramène à montrer d'une part les inclusions

$$(5.4.2) \qquad \left(z^\alpha \log^k z \,\overline{\mathbb{Q}}\{z\}^A_{-1}\right)^+ \subset \mathbb{C}\frac{1}{z^{\alpha+1}} \cdot \overline{\mathbb{Q}}\left\{\frac{1}{z}\right\}^A_0 [\log z],$$

$$\left(z^\alpha \log^k z \,\overline{\mathbb{Q}}\left\{\frac{1}{z}\right\}^A_1\right)^+ \subset \mathbb{C}\frac{1}{z^{\alpha+1}} \cdot \overline{\mathbb{Q}}\{z\}^A_0 [\log z] \ ,$$

et d'autre part que:

(5.4.3) Les seuls éléments de $\mathbb{C}[[z]][\log z]$ (resp. $\{\frac{1}{z}\}\mathbb{C}[[\frac{1}{z}]][\log z]$) dont l'image par $^+$ est dans le compositum $\mathbb{C}[\log z] \cdot \overline{\mathbb{Q}}\{z\}^A_0$ (resp. $\mathbb{C}[\log z] \cdot \overline{\mathbb{Q}}\{z\}^A_0$) appartiennent au compositum

$$\mathbb{C}[\log z] \cdot \overline{\mathbb{Q}}\{z\}^A_{-1} \quad (\text{resp. } \mathbb{C}[\log z] \cdot \overline{\mathbb{Q}}\left\{\frac{1}{z}\right\}^A_1).$$

Soit donc $\underline{a} = (a_0, a_1, a_2, \ldots)$ une suite de nombres complexes algébriques vérifiant la condition $(\mathbf{G})$, considérons les séries $F_{\underline{a}} = \sum_{n\geq 0} \frac{a_n}{n!} z^n \in \overline{\mathbb{Q}}\{z\}^A_{-1}$ et $\check{\mathfrak{f}}_{\underline{a}} = \sum_{n\geq 0} n! a_n z^{-n} \in \overline{\mathbb{Q}}\{\frac{1}{z}\}^A_1$, et posons, suivant (5.3.10):

$$\left(z^\alpha \log^k z F_{\underline{a}}\right)^+ = \sum_j \sum_{n\geq 0} b_{n,j} z^{-\alpha-1-n} \log^j z \ ,$$

$$\left(z^\alpha \log^k z \check{\mathfrak{f}}_{\underline{a}}\right)^+ = \sum_j \sum_{n\geq 0} c_{n,j} z^{-\alpha-1-n} \log^j z \ .$$

Quitte à retrancher de $z^\alpha F_{\underline{a}}$ (resp. $z^\alpha \check{\mathfrak{f}}_{\underline{a}}$) un élément de $\overline{\mathbb{Q}}[\frac{1}{z}]$ (resp. $\overline{\mathbb{Q}}[z]$), on peut supposer, pour l'étude de $(z^\alpha \log^k z \check{\mathfrak{f}}_{\underline{a}})^+$, que $\alpha \notin -\mathbb{N}$ (resp. $\alpha \notin \mathbb{N}$). Les coefficients de la puissance maximale du logarithme sont donnés par:

$$(5.4.4) \qquad b_{n,k} = (-1)^k \Gamma(\alpha+1) \frac{(\alpha+1)_n}{n!} a_n$$

$$c_{n,k} = (-1)^{k+n} \Gamma(\alpha+1) \frac{n!}{(-\alpha-1)_n} a_n \quad \text{si } \alpha \text{ n'est pas entier,}$$

$$c_{n,k+1} = (-1)^{\alpha+n} \frac{1}{(k+1)\binom{n}{-\alpha-a}} a_n \quad \text{si } \alpha \text{ est entier} < 0.$$

Comme $\alpha$ est rationnel, on en déduit déjà que les séries $\sum_{n\geq 0} b_{n,k} z^{-\alpha-1-n}$ et $\sum_{n\geq 0} c_{n,k} z^{-\alpha-1+n}$ sont dans $\mathbb{C}\frac{1}{z^{\alpha+1}} \cdot \overline{\mathbb{Q}}\{\frac{1}{z}\}^A_0$ et $\mathbb{C}\frac{1}{z^{\alpha+1}} \cdot \overline{\mathbb{Q}}\{z\}^{\overline{A}}_0$ respectivement.

Ecrivons encore (5.3.10) sous la forme

(5.4.5)

$$\left(z^{\alpha+m}\log^k z\right)^+ = \Gamma(\alpha+m+1)z^{-\alpha-1-m}\sum_{j=0}^{k}\rho_{m,j}\log^j z,$$

$$\rho_{m,j}\in\mathbb{C},\ \ \rho_{m,k+1}=(-1)^k$$

si $\alpha+m$ n'est pas un entier strictement négatif,

$$= \frac{(-1)^m}{(-\alpha-1-m)!}z^{-\alpha-1-m}\sum_{j=0}^{k+1}\rho_{m,j}\log^j z,$$

$$\rho_{m,j}\in\mathbb{C},\ \ \rho_{m,k+1}=\frac{(-1)^{\alpha+m}}{k+1},\ \text{sinon.}$$

Alors la formule $(z^{\alpha+m}\log^k z)^+ = -\frac{d}{dz}(z^{\alpha+m-1}\log^k z)^+$ (cf. (5.3.7)) donne la relation de récurrence

(5.4.6) $$\rho_{m,j}=\rho_{m-1,j}-\frac{j+1}{\alpha+m}\rho_{m-1,j+1}\qquad \text{si } m\neq -\alpha,$$

d'où l'on déduit que pour $|m|$ assez grand, $\rho_{m,j}$ s'écrit comme combinaison linéaire $\rho_{m,j} - \sum_{0\leq i\leq 1}r_{m,i}\rho_i$, où les $r_{m,i}$ sont des nombres complexes algébriques tels que les suites

$$(r_{0,i},r_{1,i},\ldots,r_{m,i},\ldots)\quad (\text{resp. } (r_{0,i},r_{-1,i},\ldots,r_{-m,i},\ldots))$$

vérifient la condition (**G**) (le point essentiel est que la croissance en $n$ du dénominateur commun des quantités rationnelles $\prod_{i=1}^{l}\frac{1}{\alpha+m_i}$, $l\leq k+1$ ($k$ fixé), $1\leq m_i\leq n$, est exponentielle). Par exemple, si $k=1$, et si $\alpha$ n'est pas entier, on a $\rho_{m,1}=-1$, $\rho_{m,0}=\frac{\Gamma'}{\Gamma}(\alpha+1)+\sum_{i=1}^{m}\frac{1}{\alpha+i}$. De ce que

(5.4.7) $$b_{nkj}=\frac{n!}{(-\alpha-1)_n}\rho_{n,j}a_n,$$

$$c_{n,j}=(-1)^n\frac{n!}{(-\alpha-1)_n}\rho_{-n,j}a_n \text{ si } \alpha \text{ n'est pas entier,}$$

$$c_{n,j}=\frac{1}{\binom{n}{-\alpha-1}}\rho_{-n,j}a_n \text{ si } \alpha \text{ est entier} < 0,$$

on en conclut que les séries $\sum_{n\geq 0}b_{n,j}z^{-\alpha-1-n}$ et $\sum_{n\geq 0}c_{n,j}z^{-\alpha-1+n}$ sont dans $\mathbb{C}\frac{1}{z^{\alpha+1}}\cdot\overline{\mathbb{Q}}\{\frac{1}{z}\}_0^A$ et $\mathbb{C}\frac{1}{z^{\alpha+1}}\cdot\overline{\mathbb{Q}}\{z\}_0^A$ respectivement, ce qui établit (5.4.2).

Pour prouver (5.4.3), considérons un ensemble fini de suites

$$(a_{0,k},a_{1,k},\ldots,a_{m,k},\ldots)$$

de nombres complexes (plus nécessairement algébriques), puis formons les séries

$$F_k=\sum_{n\geq 0}\frac{a_{n,k}}{n!}z^n\quad\text{et}\quad \bar{\mathsf{f}}_k=\sum_{n\geq 0}n!a_{n,k}z^{-n}.$$

Il s'agit de voir que si $\sum_k(\log^k zF_k)^+$ (resp. $\sum_k(\frac{1}{2}\log^k z\bar{f}_k)^+$) appartient à $\mathbb{C}[\log z]\cdot\overline{\mathbb{Q}}\{\frac{1}{z}\}^A_0$ (resp. à $\mathbb{C}[\log z]\cdot\overline{\mathbb{Q}}\{z\}^A_0$), alors chaque $F_k$ (resp. $\bar{f}_k$) appartient à $\mathbb{C}\cdot\overline{\mathbb{Q}}\{z\}^A_{-1}$ (resp. à $\mathbb{C}\cdot\overline{\mathbb{Q}}\{\frac{1}{z}\}^A_1$).

Par récurrence sur le degré en $\log z$, on se ramène à ne considérer que le maximum des $k$, et le résultat découle alors des première et troisième formules de (5.4.4).

5.5. *Solutions de* $\Phi$ *en* 0 *et* (*micro*)*solutions de* $\overline{\triangledown}\phi$ *en* $\infty$. On a établi en 5.1, 5.2 l'existence d'une base de solutions en 0 du $E$-opérateur $\Phi$ de la forme

$$(y_1(z),\ldots,y_\nu(z)) = (F_1(z),\ldots,F_\nu(z))\cdot z^{\Gamma_0}$$

où les $F_i$ sont dans $K[[z]]$, et $\Gamma_0$ désigne une matrice carrée d'ordre $\nu$ triangulaire supérieure à coefficients dans $\mathbb{Q}$. Il s'agit de voir que les $F_i$ sont des $E$-fonctions, ou, ce qui revient au même, que les $y_i$ sont dans NGA$\{z\}_{-1}$.

En utilisant les formules (5.3.7) et (5.3.8), on voit comme en (2.4) que $y_i^+$ est solution du $G$-opérateur $\frac{d^\rho}{dz^\rho}\cdot\overline{\triangledown}\Phi$ pour $\rho$ convenable (on pourrait prendre $\rho = 0$ si le $i$-ème terme diagonal de $\Gamma_0$ était non entier). D'après le théorème 3.5 (appliqué à l'infini), on a donc $y_i^+ \in$ NGA$\{\frac{1}{z}\}_0$. En vertu de 5.4.1, on en conclut que $y_i$NGA$\{z\}_{-1}$.

*Remarque.* On peut modifier légèrement l'argument de manière à ne pas supposer *a priori* que les exposants sont rationnels. On retrouve alors, par voie "formelle", que l'ensemble des classes non nulles d'exposants modulo $\mathbb{Z}$ de $\Phi$ en 0 coïncide avec l'ensemble analogue pour $\overline{\triangledown}\Phi$ en $\infty$.

5.6. *Solutions de* $\Phi$ *en* $\infty$ *et* (*micro*)*solutions de* $\overline{\triangledown}\Phi$ *en les singularités à distance finie*. Soit $K'$ l'extension finie de $K$ obtenue par adjonction des singularités à distance finie $\zeta_i$ de $\overline{\triangledown}\Phi$ et des exposants de Turrittin de $\Phi$ en $\infty$. On a établi en 5.1 et 5.2 (via Turrittin-Levelt) l'existence d'une base de solutions en $\infty$ du $E$-opérateur $\Phi$ de la forme

$$(\tilde{y}_1,\cdots,\tilde{y}_\nu) = \left(\hat{f}_1,\ldots\hat{f}_\nu\right)\cdot\left(\frac{1}{z}\right)^{\Gamma_\infty}\cdot e^{-\Delta z}$$

où les $\hat{f}_j$ sont dans $K'((\frac{1}{z}))$, $\Delta$ est la matrice diagonale ayant pour coefficients diagonaux les $\zeta_j$ (comptées avec multiplicité), et $\Gamma_\infty$ désigne une matrice sous forme de Jordan à coefficients diagonaux dans $K'$, commutant à $\Delta$. En particulier, $\Gamma_\infty$ est sous forme de blocs $\Gamma_{\infty,i}$, un pour chaque singularité finie de $\overline{\triangledown}\Phi$ (sans multiplicité).

Il s'agit de montrer d'une part que

(5.6.1) Pour tout $j$, après multiplication par une puissance de $z$, $\hat{f}_j(\frac{1}{z})$ devient une 3-fonction, ou, ce qui revient au même, que $\hat{y}_j := \tilde{y}_j e^{\zeta_j z}$ est dans NGA$\{\frac{1}{z}\}_1$,

et d'autre part que

(5.6.2) Les classes non nulles modulo $\mathbb{Z}$ des éléments diagonaux de $\Gamma_{\infty,i}$ coïncident avec classes non nulles des exposants de $\overline{\overline{\nabla}}\Phi$ en $\zeta_i$.

On commence par se ramener au cas $\zeta_i = 0$ de la manière suivante: le tordu $\Phi \otimes e^{-\zeta_j z}$ de $\Phi$ par $e^{-\zeta_j z}$ (c'est-à-dire l'opérateur obtenu en remplaçant $\frac{d}{dz}$ par $\frac{d}{dz} - \zeta_i$) annule au bloc $\Gamma_{\infty,i}$, et est encore de type $E$ car $\overline{\nabla}(\Phi \otimes e^{-\zeta_j z})$ n'est autre que le $G$-opérateur translaté de $\overline{\overline{\nabla}}\Phi$ par $\zeta_i$ (c'est-à-dire l'opérateur obtenu en remplaçant $z$ par $z + \zeta_i$). Fixant $i$, nous supposons donc que $\Gamma_{\infty,i}$ est le bloc attaché à la singularité 0 de $\overline{\overline{\nabla}}\Phi$, et nous considérons les éléments diagonaux $\alpha_j$ de $\Gamma_{\infty,i}$ et les "séries" $\hat{y}_j = \tilde{y}_j$ correspondantes, $\tilde{y}_j \in (\frac{1}{z})^{\alpha_j} K'((\frac{1}{z}))[\log z]$.

Considérons d'abord le cas des exposants de Turrittin $\alpha_j$ non entiers. Alors en utilisant les formules (5.3.7) et (5.3.8), on voit comme ci-dessus que $\tilde{y}_j^+(-z)$ est solution du $G$-opérateur $\overline{\overline{\nabla}}\Phi$ dans $z^{\alpha_j} K'((z))[\log z]$. D'après 3.5 (appliqué en l'infini), on a donc $\tilde{y}_j^+ \in \mathrm{NGA}\{z\}_0$. En vertu de 5.4.1, on en conclut que $\tilde{y}_j \in \mathrm{NGA}\{\frac{1}{z}\}_{+1}$. Réciproquement, si $y \in z^{\alpha} K'((z))[\log z]$ est solution en 0 de $\overline{\overline{\nabla}}\Phi$, avec $\alpha \notin \mathbb{Z}$, alors $y^+$ est solution dans $(\frac{1}{z})^{\alpha_j} K'((\frac{1}{z}))[\log z]$ de $\Phi$ en $\infty$. Ceci prouve (5.6.2) (en particulier $K' = K(\zeta_i)_i$), et aussi (5.6.1) dans le cas des exposants de Turrittin non entiers. La difficulté avec les exposants entiers vient bien entendu de ce qu'une série infinie apparaît au second membre de (5.3.8). Pour la contourner , il suffit en fait de multiplier $\tilde{y}_j(\frac{1}{z}) \in K'((\frac{1}{z}))[\log z]$ par $(\frac{1}{z})^{\alpha}$, pour $\alpha$ rationnel convenable, à condition de prouver que $(\frac{1}{z})^{\alpha} \tilde{y}_j$ est encore solution d'un certain $E$-opérateur. Par récurrence sur le degré du logarithme comme en 4.2, on se ramène à prouver l'énoncé suivant:

PROPOSITION 5.6.3.    Soient $\hat{f} \in \mathbb{C}((\frac{1}{z}))$ et $\alpha \in \mathbb{Q}$. Alors $\hat{f}$ une solution d'un $E$-opérateur $\Leftrightarrow z^{\alpha}\hat{f}$ est solution d'un $E$-opérateur.

On peut supposer $\hat{f} \in \frac{1}{z}\mathbb{C}[[\frac{1}{z}]]$. Prouvons $\Rightarrow$ ($\Leftarrow$ est analogue ou même plus facile: pour $\alpha \in \mathbb{Q} - \mathbb{Z}$, on pourrait l'obtenir par voie formelle en comparant $\hat{f}^+$ et $(z^{\alpha}\hat{f})^+$). Nous procéderons par voie analytique. Comme tout $E$-opérateur $\Phi$ annulant $\hat{f}$ a une seule pente non nulle à l'infini, égale à 1, on sait que $\hat{f}$ est 1-sommable au sens de Ramis dans toute direction $\theta$ non singulière, i.e. asymptotique à une (unique) fonction holomorphe $f_{\theta}$ sur un secteur $V$ d'ouverture $> \pi$ bissecté par $\theta$; en outre, $f_{\theta}$ est automatiquement solution de $\Phi$. D'après ce qui a été prouvé en 5.5, $f_{\theta}$ est donc restriction sur $V$ d'une (unique) solution $y$ de $\Phi$ dans $\mathrm{NGA}\{z\}_{-1}$. D'après 4.2, $z^{\alpha}y$ est solution d'un $E$-opérateur $\Phi'$. Son développement asymptotique dans $V$ n'est autre que $z^{\alpha}\hat{f}$, et c'est une solution de $\Phi'$ (noter que quitte à changer légèrement $\theta$, on peut supposer que $\theta$ n'est pas direction singulière de $\Phi'$). Ceci achève la preuve de 4.3.

*Remarque.* Soit $y$ une solution d'un $E$-opérateur $\Phi$ sur le demi-axe réel positif, bornée en 0. Supposons que les monodromies locales de $\overline{\nabla}\Phi$ en les singularités à distance finie $\zeta_i$ soient finies. Renumérotons $\zeta_1, \ldots, \zeta_r$ celles de ces singularités de partie réelle maximale ($r > 1$ si et seulement si l'axe réel est une ligne de Stokes).

On peut développer (le prolongement analytique de) la transformée de Laplace $y^+$ de $y$ au voisinage de $-\zeta = -\zeta_1, \ldots, -\zeta_r$ sous la forme

$$y^+(z - \zeta) \sum_{\alpha\zeta} z^{\alpha\zeta} \sum_{n\geq 0} a_{n,\alpha\zeta} z^n, \text{ pour } |z| < \varepsilon,$$

où les $\alpha_\zeta$ sont dans des classes distinctes modulo $\mathbb{Z}$. Un calcul intégral long mais élémentaire ([DiP, II8]), basé sur le lemme de Jordan, donne alors le développement asymptotique explicite de $y(z)$ pour $z \to \infty$:

$$\sum_{\zeta = \zeta_1, \ldots, \zeta_r} e^{-\zeta z} \sum_{\alpha\zeta} \left(\frac{1}{z}\right)^{\alpha\zeta} \sum_{n\geq 0} \frac{1}{\Gamma(-n - \alpha\zeta)} a_{n,\alpha\zeta} \left(\frac{1}{z}\right)^{n+1}.$$

Comme le développement de $y^+(z - \zeta)$ est élément de $\mathrm{NGA}\{z\}_0$ d'après le "théorème de permanence", on voit directement que les coefficients des $e^{-\zeta z}$ dans ce développement asymptotique sont éléments de $\mathrm{NGA}\{\frac{1}{z}\}_1$. C'est cette remarque qui m'a permis de deviner le "théorème de dualité".

Il est d'ailleurs plausible que l'on puisse tirer une preuve "analytique" de 4.3 (sans recours au calcul opérationnel) en amplifiant cette remarque.

5.7. *Preuve de 4.6.* Le cas d'un élément holonome de $\mathrm{NGA}\{z\}_1$ se ramène à celui d'une 3-fonction comme dans 4.2. Soit donc $\mathfrak{f}(z) \in K[[z]]$ une 3-fonction. Soit $\alpha \in \mathbb{Q} - \mathbb{Z}$. D'après 5.6.3 (implication $\Leftarrow$), il suffit de montrer que la série holonome $z^\alpha \hat{f}(\frac{1}{z}) \in z^\alpha K\{\frac{1}{z}\}^A_{+1}$ est solution d'un $E$-opérateur. En vertu de 5.4.1, et des formules (5.3.7) et (5.3.8) , $(z^\alpha \hat{f}(\frac{1}{z}))^+$ est un élément holonome de $z^{-\alpha-1} K\{z\}^A_0$. C'est donc une solution d'un $G$-opérateur $\Psi$ (cf. 3.2, amplifié en 3.6). Les mêmes formules (5.3.7) et (5.3.8) montrent alors que $z^\alpha \hat{f}(\frac{1}{z})$ est solution du $E$-opérateur $\overline{\nabla}\Psi$.

## 6. Pureté et dualité

THÉORÈME 6.1.     *Soient $s$ un nombre rationnel non nul, et $y(z)$ un élément holonome de $\mathrm{NGA}\{z\}_s$. Alors $y(z^{-s})$ est solution d'un $E$-opérateur différentiel. Plus précisément, si $\Psi$ est un opérateur d'ordre minimal annulant $y(z)$, il existe un $E$-opérateur différentiel qui annule toutes les solutions de $\Psi$ en 0 recalibrées par $z \to z^{-s}$.*

*Preuve.* Remarquons que pour $s < 0$ (resp. $s > 0$), $y(z^{-s})$ est un élément holonome de $\mathrm{NGA}\{z\}_{-1}$ (resp. $\mathrm{NGA}\{\frac{1}{z}\}^{+1}$); cf. 1.4.1. La première assertion est donc conséquence de la forme forte de 4.2 (resp. 4.6). Soit $\Phi$ un $E$-opérateur annulant $y(z^{-s})$.

Pour la seconde assertion, écrivons $s = p/q$ comme fraction irréductible, et considérons un opérateur $\Psi' \in \mathbb{C}[z, \frac{d}{dz}]$ d'ordre minimal annulant $y(z^{\frac{1}{q}})$. Un peu de théorie de Galois différentielle montre que $\Psi'$ annule toutes les solutions de $\Psi$ en 0 recalibrées par $z \to z^{\frac{1}{q}}$. D'autre part, le module différentiel $M_{\Psi'} := \mathbb{C}(z)[\frac{d}{dz}]/\mathbb{C}(z)[\frac{d}{dz}]\Psi'$ est quotient de $\rho_* M_\Phi$, en notant $\rho$ l'application $z \to z^{-p}$, donc $\rho^* M_{\Psi'}$ est quotient de $\rho^* \rho_* M_\Phi$. On voit par là que tout $E$-opérateur $\Phi'$ multiple commun à gauche de $E$-opérateurs annulant les $y(\varepsilon z^{-s})$, où $\varepsilon$ parcourt les racines de l'unité d'ordre $p$, annule toutes les solutions de $\Psi$ en 0 recalibrées par $z \to z^{-s}$.

6.2. Le théorème de pureté (en 0) et de dualité (pureté en $\infty$) suivent alors directement de 6.1 et 4.3 (voir aussi 4.5), sauf en ce qui concerne les assertions relatives aux développements asymptotiques.

Ces dernières découlent aisément d'une part de ce qu'il n'y a qu'une seule pente non nulle en jeu, de sorte que toutes les séries divergentes qui apparaissent sont $\frac{1}{|s|}$-sommables selon toute direction non singulière (sur de "grands" secteurs d'ouverture $> |s|\pi$; cf. [R3]), et d'autre part de ce que la $\frac{1}{|s|}$-sommation respecte la structure d'algèbre différentielle.

*Remarque.* On pourrait arithmétiser davantage notre cadre en remplaçant, dans les définitions et le théorème de dualité, le corps $\mathbb{C}$ par le sous-corps obtenu en adjoignant à $\overline{\mathbb{Q}}$ les valeurs de la fonction gamma et de ses dérivées aux points rationnels.

6.3. *Exemple*: *l'équation différentielle de Weber* [WW].

$$\frac{d^2}{dz^2}D_m(z) = \left(\frac{z^2}{4} - \frac{1}{2} - m\right)D_m(z)$$

ne présente qu'une seule singularité (à l'infini), de pente 2. Supposons $m$ rationnel mais non entier naturel, une base de solutions en 0 est alors donnée par les fonctions cylindriques paraboliques

$$D_m(\pm z) = \frac{\sqrt{\pi}2^{\frac{m}{2}}}{\Gamma\left(\frac{1-m}{2}\right)}e^{-\frac{z^2}{4}}\,_1F_1\left(\frac{-m}{2}; \frac{1}{2}; \frac{z^2}{2}\right)$$
$$-\frac{\sqrt{\pi}2^{\frac{m+1}{2}}}{\Gamma\left(\frac{-m}{2}\right)}(\pm z)e^{-\frac{z^2}{4}}\,_1F_1\left(\frac{1-m}{2}; \frac{3}{2}; \frac{z^2}{2}\right)$$

qui appartiennent à $\mathrm{NGA}\{z\}_{-\frac{1}{2}}$.

D'autre part, dans tout secteur d'ouverture $< \pi/2$ bissecté par le demi-axe réel positif, $D_m(\pm z)$ admet le développement asymptotique

$$(\pm z)^m e^{-\frac{z^2}{4}} {}_2F_0 \left( \frac{1-m}{2}, \frac{-m}{2}; -2 \left( \frac{1}{z} \right)^2 \right),$$

dont on vérifie directement qu'il satisfait l'équation de Weber dans l'algèbre différentielle $z^{m\mathbb{Z}} \exp(\frac{z^2}{4} - \mathbb{Z}) \mathbb{Q}((\frac{1}{z}))$. Pour $D_m(+z)$, ce développement asymptotique est du reste valide dans tout secteur d'ouverture $< 3\pi/2$ bissecté par le demi-axe réel positif, de sorte que la 2-sommation au sens de Ramis de

$$ {}_2F_0 \left( \frac{1-m}{2}, \frac{-m}{2}; -2 \left( \frac{1}{z} \right)^2 \right) $$

selon toute direction d'angle $\in ]-\pi/2, \pi/2[$ est $z^{-m} e^{\frac{z^2}{4}} D_m(z)$. Par contre, lorsqu'on franchit la ligne de Stokes d'angle $-\pi/4$, le développement asymptotique de $D_m(-z)$ doit être corrigé par addition du terme

$$\frac{-\sqrt{2\pi}}{\Gamma(-m)} e^{\frac{z^2}{4}} z^{-m-1} {}_2F_0 \left( 1 + \frac{m}{2}, \frac{1+m}{2}; 2 \left( \frac{1}{z} \right)^2 \right),$$

qui satisfait aussi l'équation de Weber dans l'algèbre différentielle

$$z^{m\mathbb{Z}} \exp\left( \frac{z^2}{4} \mathbb{Z} \right) \mathbb{Q} \left( \frac{1}{z} \right).$$

Remarquons pour finir que

$$ {}_2F_0 \left( \frac{1-m}{2}, \frac{-m}{2}; -2 \left( \frac{1}{z} \right)^2 \right) $$

et

$$ {}_2F_0 \left( 1 + \frac{m}{2}, \frac{1+m}{2}; 2 \left( \frac{1}{z} \right)^2 \right) $$

sont dans $\mathbb{Q}\{\frac{1}{z}\}^A_{+\frac{1}{2}}$.

## 7. Séries de factorielles et équations aux différences finies

7.1.  On note $\Delta$ l'opérateur de différences finies de pas 1 : $(\Delta g)(x) = g(x+1) - g(x)$. Considérons l'espace vectoriel complexe $\mathbb{C}[!x!]^{(\rho)}$ des "*séries de factorielles* (inverses) *généralisées*"

$$g(x) = \sum_{n \geq 0} b_n \frac{\Gamma(x)}{\Gamma(x+n+\rho+1)},$$

où l'"exposant" $\rho$ est un nombre complexe donné. Cet espace est en fait un $\mathbb{C}[x, \Delta]$-module. Le cas des séries de factorielles inverses, introduites par Nicole dans son "traité du calcul des différences finies" (Mém. Acad. Sci. Paris 1717), correspond à $\rho = 0$, i.e. $g(x) = \sum_{n \geq 0} \frac{b_n}{x(x+1)\cdots(x+n)}$. Rappelons que dans ce cas, le développement formel à l'infini $g(x) = \sum_{n \geq 0} a_n x^{-n-1}$ peut se calculer ainsi: en posant $F(z) = \sum_{n \geq 0} \frac{a_n}{n!} z^n$ (de sorte que $\sum_{n \geq 0} a_n z^{-n-1} = F^+$), on a $\sum_{n \geq 0} \frac{b_n}{n!} z^n = F(\log \frac{1}{1-z})$ (cf. [WW, 7.82]).

La *transformée de Mellin* formelle est l'application linéaire

$$\mathcal{M} : \mathbb{C}[!x!]^{(\rho)} \to (1-z)^\rho \mathbb{C}[[1-z]]$$

définie par

$$\mathcal{M}(g) = (1-z)^\rho \sum_{n \geq 0} b_n \frac{(1-z)^n}{\Gamma(n+\rho+1)}$$

(cf. [BaD]). Cette application est injective si $\rho$ n'est pas un entier strictement négatif. Un élément de $\mathbb{C}[!x!]^{(\rho)}$ est dit Gevrey d'ordre $s$ si la série $\sum_{n \geq 0} b_n \frac{(1-z)^n}{\Gamma(n+\rho+1)}$ est Gevrey d'ordre $s$ au sens usuel (en la variable $1 - z$). D'autre part, la transformée de Mellin opératorielle est l'application $\mathcal{M} : \mathbb{C}[x, \Delta] \to \mathbb{C}[z, \frac{d}{dz}]$ définie par $x \to -z\frac{d}{dz}$, $\Delta \to z - 1$. C'est un isomorphisme de $\mathbb{C}$-algèbres; en outre, si $g$ est annulé par un élément $\Xi$ de $\mathbb{C}[x, \Delta]$, alors $\mathcal{M}(g)$ est annulé par $\mathcal{M}(\Xi)$ (loc. cit.).

7.2. Une série de factorielles généralisée g(x) est dite *Gevrey d'ordre s de type arithmétique* (pour des nombres rationnel $s$ et $\rho$ fixés) si les $b_n$ sont algébriques et la suite de terme général $\frac{b_n}{(n!)^{s+1}}$ vérifie la condition (**G**) de 1.1. (en particulier $g(x)$ est Gevrey d'ordre $s$ au sens ci-dessus). Pour $\rho = 0$, une telle série n'est donc autre qu'un élément de $\mathcal{M}^{-1}(\overline{\mathbb{Q}}\{1-z\}_s^A)$.

Du théorème de pureté pour les équations différentielles, on déduit alors, via la transformation de Mellin, le résultat suivant.

THÉORÈME 7.3. *Soit $\Xi$ un élément non nul de $\mathbb{C}[x, \Delta]$, et soient $g$ et $g'$ des solutions de $\Xi$ dans $\mathbb{C}[!x!]^{(\rho)}$ et $\mathbb{C}[!x!]^{(\rho')}$ respectivement, avec $\rho, \rho' \in \mathbb{C} \backslash \mathbb{Z}_{<0}$. On suppose que $\Xi$ est d'ordre minimal en $x$ parmi les opérateurs aux différences finies polynômiaux annulant $g$. Alors si $g$ est Gevrey d'ordre $s$ de type arithmétique, il en est de même de $g'$.*

Noter que la rationalité de $\rho$ fait partie de l'hypothèse, et celle de $\rho'$ de la conclusion.

7.4. *Remarques et perspectives.* Je considère les résultats de cette étude comme la pointe d'un iceberg de questions inexplorées, dont les plus patentes sont:

a) Généralisation en dimension supérieure, en vue d'englober par exemple des fonctions de type hypergéométrique beaucoup plus générales. Déjà pour $s = 0$, les résultats du paragraphe 3 ne sont pas complètement écrits à plusieurs variables.[13]

b) Peut-on définir une "bonne classe" d'opérateurs différentiels de type arithmétique, en étudiant l'indice à l'origine et à l'infini dans les espaces de séries Gevrey de type arithmétique?

c) Est-il nécessaire de se limiter à des équations différentielles linéaires? On sait en effet que la théorie Gevrey "complexe" a d'importantes applications non-linéaires; voir [R3] pour un tour d'horizon, et les travaux d'Ecalle.

d) Y a-t-il des $q$-analogues? Par exemple, peut-on considérer la fonction de Tschakaloff [T] $\sum_{n \geq 0} q^{-n(n-1)/2} z^n$ comme une $q$-$E$-fonction? Rappelons qu'il y a une théorie $q$-Gevrey "complexe" [Bé], [R3], [MaZ], et en particulier une transformation $q$-Laplace (et même plusieurs); nous reviendrons en partie sur cette question à la fin du second volet de cet article.[14]

e) Combiner c) et d) serait du reste intéressant au vu des nouvelles fonctions introduites par L. Denis [De], qui "interpolent" l'exponentielle de Drinfeld et l'exponentielle usuelle, et satisfont des équations aux $q$-différences polynômiales.

f) La théorie des équations différentielles linéaires $p$-adiques présente des phénomènes de monodromie "très sauvage" liés aux exposants de Liouville (sans doute absents des équations définies sur $\overline{\mathbb{Q}}$), mais aussi des phénomènes de monodromie "sauvage" plus subtils , étudiés par Christol et Mebkhout, qui ont lieu même pour des équations définies sur $\overline{\mathbb{Q}}$ à exposants de Turrittin rationnels (par exemple des modules exponentiels de Dwork), où il arrive que le rayon de convergence des développements de Turrittin soit inférieur au rayon "attendu". Il est plausible que pour les équations différentielles "de type arithmétique" considérées dans cet article, ces phénomènes n'apparaissent que pour un nombre fini de premiers $p$ (cf. 4.8).

Institut de Mathématiques, Paris, France
*E-mail address*: andre@math.jussieu.fr

---

[13] Depuis la soumission de cet article, L. Di Vizio a généralisé en dimension supérieure les théorèmes 3.1 à 3.3 ci-dessus.

[14] Sur cette question aussi, de récents progrès significatifs ont été accomplis par L. Di Vizio.